\documentclass{amsart}

\DeclareMathSymbol{\twoheadrightarrow}
{\mathrel}{AMSa}{"10}

\def\Q{{\mathbf Q}}
\def\Z{{\mathbf Z}}
\def\C{{\mathbf C}}

\def\F{{\mathbf F}}
\def\L{{\mathbf L}}
\def\H{{\mathrm H}}
\def\SS{{\mathbf S}}
\def\A{{\mathbf A}}
\def\Sn{{\mathbf S}_n}
\def\An{{\mathbf A}_n}
\def\alg{\mathrm{alg}}
\def\Gal{\mathrm{Gal}}
\def\Perm{\mathrm{Perm}}

\def\Lie{\mathrm{Lie}}
\def\Ad{\mathrm{Ad}}
\def\MT{\mathrm{MT}}
\def\mt{\mathrm{mt}}

\def\End{\mathrm{End}}
\def\Aut{\mathrm{Aut}}

    \def\sll{\mathfrak{sl}}
\def\I{\mathrm{Id}}

\def\s{{\mathcal S}}
\def\fchar{\mathrm{char}}

\def\GL{\mathrm{GL}}
\def\PGL{\mathrm{PGL}}
\def\PSL{\mathrm{PSL}}
\def\SL{\mathrm{SL}}
\def\Sp{\mathrm{Sp}}
\def\Gp{\mathrm{Gp}}
    \def\sp{\mathfrak{sp}}
\def\M{\mathrm{M}}
\def\E{\mathrm{E}}
\def\dim{\mathrm{dim}}

\def\P{{\mathbf P}}

\def\g{{\mathfrak g}}
\def\c{{\mathfrak c}}
\def\s{{\mathfrak s}}

\def\G{{\mathfrak G}}
\def\K{{\mathfrak K}}

\def\m{{\mathfrak m}}

\def\ZZ{{\mathfrak Z}}

\newtheorem{thm}{Theorem}[section]
\newtheorem{lem}[thm]{Lemma}
\newtheorem{cor}[thm]{Corollary}
\newtheorem{prop}[thm]{Proposition}
\theoremstyle{definition}
\newtheorem{defn}[thm]{Definition}

\newtheorem{rem}[thm]{Remark}
\newtheorem{rems}[thm]{Remarks}
\title[Hyperelliptic jacobians]
{Very simple $2$-adic representations and hyperelliptic jacobians}
\author[Yuri G. Zarhin]{Yuri G. Zarhin}
\thanks{Partially supported by the NSF}

\begin{document}
\maketitle
\section{Introduction}
Throughout this paper $K$ is a field of characteristic $\ne 2$ and
$K_a$ its algebraic closure. If $f(x) \in K[x]$ is a separable
polynomial of degree $n \ge 5$ then it gives rise to the
hyperelliptic curve
$$C=C_f:y^2=f(x).$$
We write $J(C)=J(C_f)$ for its jacobian; it is a $g$-dimensional
abelian variety defined over $K$, whose dimension $g$ is equal to
$\frac{n-1}{2}$ when $n$ is odd and equal to $\frac{n-2}{2}$ when
$n$ is even. In this paper we deal with special (but somehow
``generic" ) case when the Galois group $\Gal(f)$ of  $f$ is
either the symmetric group $\Sn$ or the alternating group $\An$.
We study the $\ell$-adic Lie algebra $g_{\ell,J(C)}$ attached to
the Galois action on the $\ell$-torsion of $J(C)$ and prove that
this Lie algebra is ``as large as possible" when $K$ is a number
field. As a corollary, we obtain the Tate conjecture  and the
Hodge conjecture for all self-products of $J(C)$. In fact, we
prove that all the Tate/Hodge classes involved can be presented as
linear combinations of products of divisor classes.

Our approach is based on a study of the $2$-adic image of the
Galois group $\Gal(K)$ in the automorphism group
$\Aut_{\Z_2}(T_2(J(C)))$ of the $\Z_2$-Tate module $T_2(J(C))$ of
$J(C)$. We prove that the {\sl algebraic envelope} of the image
contains  the symplectic group attached to the theta divisor on
$J(C_f)$. Our proof is based on known lower bounds for dimensions
of nontrivial (projective) representations of $\A_n$ in
characteristics $0$ and $2$ (\cite{JK}, \cite{Wagner},
\cite{Wagner2}) and a notion of {\sl very simple} representation
introduced and studied in \cite{ZarhinTexel} and
\cite{ZarhinCrelle}. This allows us to prove that $g_{\ell,J(C)}$
is ``as large as possible" for $\ell=2$. Now the rank independence
on $\ell$  for $g_{\ell,J(C)}$ \cite{ZarhinInv} allows us to to
extend this assrtion to all primes $\ell$, using a variant of a
theorem of Borel - de Siebenthal \cite{ZarhinDuke}.

The present paper is a natural continuation of our previous
articles \cite{ZarhinIzv}, \cite{ZarhinMRL}, \cite{ZarhinTexel},
\cite{ZarhinMRL2} and \cite{ZarhinCrelle}.

\section{Abelian varieties and $\ell$-adic Lie algebras}
\label{mainr} Let $\ell$ be a prime and $K$ be a field of
characteristic different from $\ell$. We fix its algebraic closure
$K_a$ and write $\Gal(K)$ for the absolute Galois group
$\Aut(K_a/K)$. Let $K({\ell})$ be the abelian extension of $K$
obtained by adjoining to $K$ all ${\ell}$-power roots of unity. We
write
$$\chi_{\ell}: \Gal(K) \twoheadrightarrow \Gal(K({\ell})/K) \subset \Z_{\ell}^*$$
for the corresponding cyclotomic character. If every finite
algebraic extension of $K$ contains only finitely many
$\ell$-power roots of unity (e.g., $K$ is finitely generated over
its prime subfield) then the image $\chi_{\ell}(\Gal(K))$ is
infinite. We write $\Z_{\ell}(1)$ for the $\Gal(K)$-module
$\Z_{\ell}$ with the Galois action defined by character
$\chi_{\ell}$. If $X$ is an abelian variety defined over $K$ and
$m$ is a positive integer not divisible by $\fchar(K)$ then we
write $X_m$ for the kernel of multiplication by $m$ in $X(K_a)$.
It is well-known \cite{MumfordAV} that $X_m$ is a free
$\Z/m\Z$-module of rank $2\dim(X)$ provided with a natural
structure of $\Gal(K)$-module. Suppose $\ell$ is a prime distinct
from $\fchar(K)$. As usual, we write $T_{\ell}(X)$ for the
projective limit of $\Gal(K)$-modules $X_{\ell^i}$ where the
transition maps are multiplications by $\ell$. It is well-known
that $T_{\ell}(X)$ is a free $\Z_{\ell}$-module of rank $2\dim(X)$
provided with natural continuous homomorphism ($\ell$-adic
representation)
$$\rho_{\ell,X}: \Gal(K) \to \Aut_{\Z_{\ell}}(T_{\ell}(X)) \cong
\GL(2\dim(X),\Z_{\ell}).$$ In addition $X_{\ell}=T_{\ell}(X)/\ell
T_{\ell}(X)$ (as Galois module).

 Recall \cite{MumfordAV} that
each   polarization $\lambda$ on $X$ defined over $K$ gives rise
to Riemann form

$$e_{\lambda}: T_{\ell}(X) \times T_{\ell}(X) \to \Z_{\ell}(1)\cong \Z_{\ell};$$
$e_{\lambda}$ is a non-degenerate $\Gal(K)$-equivariant
alternating $\Z_{\ell}$-bilinear form on $T_{\ell}(X)$. It is
perfect if and only if $\ell$ does not divide $\deg(\lambda)$
(e.g., $\lambda$ is principal).
 Clearly,
the corresponding automorphism group
$$\Aut_{\Z_{\ell}}(T_{\ell}(X), e_{\lambda})\cong
\Sp(2\dim(X),\Z_{\ell}).$$

It is well-known that the image
$$G_{{\ell},X}:=\rho_{{\ell},X}(\Gal(K))\subset \Aut_{\Z_{\ell}}(T_{\ell}(X))$$
is contained in the group of {\sl symplectic similitudes}
$$\Gp(T_{\ell}(X),e_{\lambda}):=\{u \in\Aut_{\Z_{\ell}}(T_{\ell}(X))\mid \exists
c=c_u \in \Z_{\ell}^* \text{ such that }$$
$$e_{\lambda}(ux,uy)=c
u(x,y) \quad \forall x,y \in T_{\ell}(X)\}.$$ More precisely, if
$\sigma \in \Gal(K)$ and $u=\rho_{{\ell},X}(\sigma)$ then
$c_u=\chi_{\ell}(\sigma)$.
It is also well-known that the
composition of $\rho_{{\ell},X}$ and the determinant map
$$\det: \Aut_{\Z_{\ell}}(T_{\ell}(X)) \to \Z_{\ell}^*$$
coincides with $\chi_{\ell}^{\dim(X)}$.

 We write $\tilde{\rho}_{\ell,X}$ for the corresponding modular
representation
$$\tilde{\rho}_{\ell,X}: \Gal(K) \to \Aut(X_{\ell});$$
we denote by $\tilde{G}_{\ell,X}$ the image
$\tilde{\rho}_{\ell,X}(\Gal(K)) \subset \Aut(X_{\ell})$. Clearly,
$\tilde{\rho}_{\ell,X}$ coincides with the composition of
${\rho}_{\ell,X}$ and the reduction map modulo $\ell$
$$\Aut_{\Z_{\ell}}(T_{\ell}(X))\twoheadrightarrow \Aut(X_{\ell});$$
 the subgroup $\tilde{G}_{\ell,X}$ coincides with the image of
${G}_{\ell,X}$ under the reduction map.

 As usual, $V_{\ell}(X)$ stands for the $\Q_{\ell}$-Tate module
$$V_{\ell}(X)=T_{\ell}(X)\otimes_{\Z_{\ell}}\Q_{\ell};$$
it is a $2\dim(X)$-dimensional $\Q_{\ell}$-vector space provided
with a natural structure of $\Gal(K)$-module and $T_{\ell}(X)$ is
identified with a $\Gal(K)$-stable $\Z_{\ell}$-lattice in
$V_{\ell}(X)$. The form $e_{\lambda}$ extends by
$\Q_{\ell}$-linearity to the non-degenerate alternating
$\Q_{\ell}$-bilinear form
$$e_{\lambda}:V_{\ell}(X)\times V_{\ell}(X) \to \Q_{\ell}.$$
We write $\Sp(V_{\ell}(X),e_{\lambda})$ for the corresponding
symplectic group viewed as $\Q_{\ell}$-linear algebraic subgroup
of $\GL(V_{\ell}(X)$. Let
$$\sp(V_{\ell}(X),e_{\lambda}) \cong \sp(2\dim(X), \Q_{\ell})$$
be the Lie algebra of $\Sp(V_{\ell}(X),e_{\lambda})$; it is an
absolutely irreducible $\Q_{\ell}$-linear subalgebra of
$\End_{\Q_{\ell}}(V_{\ell}(X))$.

It is well-known \cite{SerreLocal} that $G_{{\ell},X}$ is an
${\ell}$-adic Lie subgroup of $\Aut_{\Z_{\ell}}(T_{\ell}(X))$. We
write $g_{{\ell},X}$ for the Lie algebra of of $G_{{\ell},X}$; it
is a $\Q_{\ell}$-Lie subalgebra of
$\End_{\Q_{\ell}}(V_{\ell}(X))$. The inclusion
$G_{{\ell},X}\subset \Gp(T_{\ell}(X),e_{\lambda})$ implies easily
that
$$g_{{\ell},X}\subset \Q_{\ell} \I \oplus \sp(V_{\ell}(X),e_{\lambda}).$$
In addition, $g_{{\ell},X}\subset  \sp(V_{\ell}(X),e_{\lambda})$
if either $K=K({\ell})$ or there exists a finite algebraic
extension of $K$ that contains all $\ell$-power roots of unity.
Clearly, $g_{{\ell},X}=\{0\}$ if and only if $G_{{\ell},X}$ is
finite. It is also clear that if  every finite algebraic extension
of $K$ contains only finitely many $\ell$-power roots of unity
then $g_{{\ell},X}$ does {\sl not} lie in the Lie algebra
$\sll(V_{\ell}(X))$ of linear operators in $V_{\ell}(X)$ with zero
trace.
\begin{rem}
\label{bog} If $K$ is a field of characteristic zero finitely
generated over $\Q$ then, by a theorem of Bogomolov
\cite{Bogomolov}, $g_{\ell,X}$ contains $\Q_{\ell} \I$.
\end{rem}

\begin{thm}
\label{main0} Let us assume that $\dim(X) \ge 4$ and let us put
$d=2\dim(X)$. Let us assume that $\fchar(K)\ne 2$ and let us put
$\ell=2$. Suppose  $\tilde{G}_{2,X}$ contains a subgroup
isomorphic to the alternating group $\A_{d+1}$ (e.g.,
$\tilde{G}_{2,X}=\A_{d+1}, \SS_{d+1}, \SS_{d+2}$ or $\A_{d+2}$).
Then either $g_{2,X}=\Q_{2}\I\oplus\sp(V_{2}(X),e_{\lambda})$ or
$g_{2,X}=\sp(V_{2}(X),e_{\lambda})$. In addition, the ring
$\End(X)$ of all $K_a$-endomorphisms of $X$ is $\Z$.

If every finite algebraic extension of $K$ contains only finitely
many $2$-power roots of unity then  $g_{2,X}=\Q_{2}\I\oplus
\sp(V_{2}(X),e_{\lambda})$. If there exists a finite algebraic
extension of $K$ that contains all $\ell$-power roots of unity
then $g_{2,X}=\sp(V_{2}(X),e_{\lambda})$.
\end{thm}

\begin{thm}
\label{main1} Suppose $K$ is a field and $\fchar(K)\ne 2$. Suppose
$X$ is an abelian variety defined over $K$.
 Let us assume that $\dim(X) \ge 4$ and let us put
$d=2\dim(X)$.  Let us assume that $\tilde{G}_{2,X}$ contains a
subgroup isomorphic to the alternating group $\A_{d+1}$ and $K$
enjoys one of the two following properties:
\begin{enumerate}
\item[(i)]
$K$ is a field of characteristic zero finitely generated over $\Q$;
\item[(ii)]
$p=\fchar(K)>0$ and $K$ is is a global field.
\end{enumerate}
Then for all primes $\ell \ne \fchar(K)$
$$g_{{\ell},X}=\Q_{\ell}\I\oplus \sp(V_{\ell}(X),e_{\lambda}).$$
\end{thm}

\begin{thm}
\label{main} Let $K$ be a field with $\fchar(K) \ne 2$,
 $K_a$ its algebraic closure,
$f(x) \in K[x]$ a separable polynomial of  degree $n \ge 9$, whose
Galois group $\Gal(f)$ enjoys one of the following properties:
\begin{enumerate}
\item[(i)]
$\Gal(f)$ is either $\Sn$ or $\An$;
\item[(ii)]
$n=11$ and $\Gal(f)$ is the Mathieu group $\M_{11}$;
\item[(iii)]
$n=12$ and $\Gal(f)$ is either the Mathieu group $\M_{12}$ or
$\M_{11}$.
\item[(iv)]
There exist an odd power prime $q$ and an integer $m \ge 3$ such
that $(q,m) \ne (3,4)$,
 $n=\frac{q^m-1}{q-1}$ and $\Gal(f)$
contains a subgroup isomorphic  to the projective special linear
group $\L_m(q):=\PSL_m(\F_q)$. (E.g., $\Gal(f)$ is isomorphic
either to the projective linear group $\PGL_m(\F_q)$ or to
$\L_m(q)$.)
\end{enumerate}
 Let $C_f$ be the hyperelliptic
curve $y^2=f(x)$. Let $X=J(C_f)$ be its jacobian, $\lambda$ the
principal polarization on $J(C_f)$ attached to the theta divisor.
Then $\rho_{2,J(C_f)}(\Gal(K))$ contains an open subgroup of
$\Aut_{\Z_2}(T_2(X), e_{\lambda})$. In addition, the ring
$\End(X)$ of all $K_a$-endomorphisms of $X$ is $\Z$.
\end{thm}

\begin{thm}
\label{jac} Let $f(x) \in K[x]$ be a separable polynomial of
degree $n \ge 9$, whose  Galois group $\Gal(f)$ enjoys one of the
following properties:
\begin{enumerate}
\item[(i)]
$\Gal(f)$ is either $\Sn$ or $\An$;
\item[(ii)]
$n=11$ and $\Gal(f)$ is the Mathieu group $\M_{11}$;
\item[(iii)]
$n=12$ and $\Gal(f)$ is either the Mathieu group $\M_{12}$ or
$\M_{11}$.
\item[(iv)]
There exist an odd power prime $q$ and an integer $m \ge 3$ such
that $(q,m) \ne (3,4)$,
 $n=\frac{q^m-1}{q-1}$ and $\Gal(f)$
contains a subgroup isomorphic  to  $\L_m(q)$. (E.g., $\Gal(f)$ is
isomorphic either to the projective linear group $\PGL_m(\F_q)$ or
to $\L_m(q)$.)
\end{enumerate}
 Let $C_f$ be the
hyperelliptic curve $y^2=f(x)$. Let $X=J(C_f)$ be its jacobian.
Assume that either $K$ is a field of characteristic zero finitely
generated over $\Q$ or $K$ is a global field of odd
characteristic. Then for all primes $\ell \ne \fchar(K)$
$$g_{{\ell},X}=\Q_{\ell}\I\oplus \sp(V_{\ell}(X),e_{\lambda}).$$
\end{thm}

\begin{thm}
\label{jacall} Let $K$ be a  field  of characteristic zero
finitely generated over $\Q$, $f(x) \in K[x]$  a separable
polynomial of degree $n \ge 5$, whose Galois group $\Gal(f)$
enjoys one of the following properties:
\begin{enumerate}
\item[(i)]
$\Gal(f)$ is either $\Sn$ or $\An$;
\item[(ii)]
$n=11$ and $\Gal(f)$ is the Mathieu group $\M_{11}$;
\item[(iii)]
$n=12$ and $\Gal(f)$ is either the Mathieu group $\M_{12}$ or
$\M_{11}$.
\item[(iv)]
There exist an odd power prime $q$ and an integer $m \ge 3$ such
that $(q,m) \ne (3,4)$,
 $n=\frac{q^m-1}{q-1}$ and $\Gal(f)$
contains a subgroup isomorphic  to $\L_m(q)$. (E.g., $\Gal(f)$ is
isomorphic either to the projective linear group $\PGL_m(\F_q)$ or
to $\L_m(q)$.)
\end{enumerate}
 Let $C_f$ be the hyperelliptic curve $y^2=f(x)$. Let $X=J(C_f)$
be its jacobian. Then for all primes $\ell$
$$g_{{\ell},X}=\Q_{\ell}\I\oplus \sp(V_{\ell}(X),e_{\lambda}).$$
\end{thm}

We prove Theorems \ref{main0}, \ref{main1}, \ref{jac} and
\ref{jacall} in Section \ref{mainproof}.

\section{Group theory}
 Throughout the paper we will freely use the following
elementary observation.

\begin{prop}
\label{cont}
 Suppose $\ell$ is a prime, $F$ is field which is a finite algebraic extension of
 $\Q_{\ell}$. Suppose $W$ is a finite-dimensional $F$-vector space
 and $G$ is a compact subgroup of $\Aut_{F}(W)$ (in $\ell$-adic topology).
\begin{enumerate}
\item[(i)]
 Suppose $M$ is a periodic group of finite
 exponent provided with discrete topology (e.g., a finite group).
 Then every homomorphism  $\pi:G \to M$  is continuous and
 its kernel is an open subgroup of finite index in $G$.
\item[(ii)]
Suppose $H$ is an open normal subgroup of finite index in $G$ and
the quotient $G/H$ is a finite simple non-abelian group. Then:
\begin{enumerate}
\item[(a)]
$H':=H\bigcap \ker(\pi)$ is an open normal subgroup of finite
index in $G':=\ker(\pi)$. If $H$ contains $\ker(\pi)$ then $G/H$
is a homomorphic image of $\pi(M)$, i.e., there is a surjective
continuous homomorphism $\pi(M)\to G/H$. If $H$ does not contain
$\ker(\pi)$ then $H':=H\bigcap \ker(\pi)$ is an open subgroup of
finite index in $G':=\ker(\pi)$ and $G'/H'=G/H$. In other words,
$G/H$ is a homomorphic image either of $\ker(\pi)$ or of the image
of $\pi$.
\item[(b)]
If $M$ is solvable then $H$ does not contain $\ker(\pi)$. In
particular, $G'/H'=G/H$.
\item[(c)]
Suppose that either $\ell=2$ or $G/H$ is either a simple group of
Lie type in odd characteristic or one of $26$ simple sporadic
groups. (In other words, modulo the classification, if $\ell$ is
odd then $G/H$ is not a simple group of Lie type in characteristic
$2$.). Suppose also that $M$ is finite and for  some positive
integer $m$ there exists an embedding
$$M \hookrightarrow \PGL_m(\bar{\F}_{\ell}).$$
Then  either  $G'/H'=G/H$ or there exists an embedding
$$G/H \hookrightarrow \PGL_m(\bar{\F}_{\ell}).$$
\item[(d)]
Suppose $\pi':G \to M'$ is a continuous group homomorphism from
$G$ to an $\ell$-adic Lie group $M'$. Then :
\begin{enumerate}
\item[(d1)]
 $\pi'(H)$ is an open subgroup of finite index in the compact  group
$\pi'(G)$ and $H':=H\bigcap \ker(\pi')$ is an open subgroup of
finite index in the compact group $G':=\ker(\pi')$. In addition,
the (super)orders of $\pi'(H)$ and $H'$ both divide the
(super)order of $H$;
\item[(d2)]
Either $H$ contains $\ker(\pi')$ and  $\pi'(G)/\pi'(H)=G/H$ or
 $H$ does not contains $\ker(\pi')$ and $G'/H'=G/H$. In other words,
$G/H$ is a homomorphic image either of $\ker(\pi)$ or of the image
of $\pi$.
\item[(d3)]
Suppose $G$ is a closed subgroup of a product $S_1 \times \cdots
\times S_m$ of finitely many $\ell$-adic Lie groups $S_i$. Then
there exist a factor $S_j$, a compact subgroup $\bar{G} \subset
R_j$ and an open normal subgroup $\bar{H}$ of $\bar{G}$ such that
the quotient $\bar{G}/\bar{H}$ is isomorphic to $G/H$. In
addition,  one may choose $\bar{H}$ in such a way that the
(super)order of $\bar{H}$ divides the (super)order of $H$.
\end{enumerate}
\item[(e)]
Suppose $E$ is field which is a finite algebraic extension of
 $\Q_{\ell}$. Suppose $\alpha: \K_1 \to \K_2$ is a central isogeny
 between two semisimple $E$-algebraic groups $\K_1$ and $\K_2$. Suppose
 $\gamma: G \to \K_2(E)$ is a continuous group homomorphism (in
 $\ell$-adic topology), whose kernel is a finite commutative
 group. Let us put
 $$G_{\alpha}:=\gamma^{-1}(\alpha(\K_1(E)))=\{g\in G\mid
 \gamma(g)\in \alpha(\K_1(E))\subset \K_2(E)\}, \quad H_{\alpha}:=H\bigcap
 G_{\alpha};$$
$$G_{-1,\alpha}:=\{u \in \K_1(E))\mid \alpha(u) \in
\gamma(G_{\alpha})\subset \K_2(E)\},$$
$$H_{-1,\alpha}:=\{u \in
\K_1(E)\mid \alpha(u) \in \gamma(H_{\alpha})\subset \K_2(E)\}.$$
 Then:
\begin{enumerate}
\item[(e1)]
$G_{\alpha}$ is an open normal subgroup of finite index in $G$,
$H_{\alpha}$ is an open normal subgroup of finite index in
$G_{\alpha}$ and $G_{\alpha}/H_{\alpha}=G/H$.
\item[(e2)]
$G_{-1,\alpha}$ is a compact group, $H_{-1,\alpha}$ is an open
subgroup of finite index in $G_{-1,\alpha}$ and
$G_{-1,\alpha}/H_{-1,\alpha}=G/H$.
\item[(e3)]
Assume  that $H$ is a pro-$\ell$-group. Then  the (super)orders of
$H_{\alpha}$ and $\gamma(H_{\alpha})$ divide the (super)order of
$H$.  If $q$ is a prime divisor of the (super)order of
$H_{-1,\alpha}$ then either $q$ divides the (super)order of $H$ or
$q$ divides $\deg(\alpha)$.
\end{enumerate}
\end{enumerate}
\end{enumerate}
\end{prop}

\begin{proof}
Let us do the case (i). It suffices to check that $\ker(\pi)$ is
an open subgroup and therefore is closed. If $n$ is the exponent
of $M$ then $\ker(\pi)$ contains $G^n:=\{x^n\mid x \in G\}$. Let
us consider $W$ as finite-dimensional $\Q_{\ell}$-vector space.
Then $G$ is a compact $\ell$-adic Lie group, thanks to $\ell$-adic
``Cartan Theorem" (\cite{SerreLie}, Part II, Sect. 9). Now basic
properties of the exponential map imply that $G^n$ contains an
open neighborhood of the identity. This implies that $G^n$
contains an open subgroup and therefore $\ker(\pi)$ also contains
an open subgroup. Since every subgroup containing an open subgroup
is also open, $\ker(\pi)$ is an open subgroup of $G$. The
compactness of $G$ implies easily that the index is finite.

Now let us do the case (ii).
\begin{enumerate}
\item[(a)]
If $\ker(\pi) \subset H$ then we have
$$\pi(G)\cong G/\ker(\pi) \twoheadrightarrow G/H.$$
If $H$ does not contain $\ker(\pi)$ then $G'/H'$ is a {\sl
nontrivial} normal subgroup of $G/H$. Now the simplicity of $G/H$
implies that $G'/H'=G/H$.
\item[(b)]
Clearly, $\pi(G) \subset M$ is solvable. If $H$ contains
$\ker(\pi)$ then the simple non-abelian $G/H$ becomes a
homomorphic image of solvable $\pi(G)$. Contradiction.
\item[(c)]
Replacing $M$ by $\pi(G)$, we may assume that $\pi$ is surjective.
We may also assume that $H$ contains $\ker(\pi)$ and therefore
there is a surjective homomorphism $\alpha: M \twoheadrightarrow
G/H$. Since $M$ is isomorphic to a subgroup of
$\PGL_m(\bar{\F}_{\ell})$, it follows from a theorem of Feit-Tits
(\cite{FT}; see also \cite{KL}) that $G/H$ is also isomorphic to a
subgroup of $\PGL_m(\bar{\F}_{\ell})$.
\item[(d)]
The assertion (d1) is obvious. If $H$ contains $\ker(\pi')$ then
we have
$$G/H=(G/\ker(\pi'))/(H/\ker(\pi'))=\pi'(G)/\pi'(H).$$
If $H$ does not contain $\ker(\pi')$ then $G'/H'$ is a {\sl
nontrivial} normal subgroup of $G/H$. Now the simplicity of $G/H$
implies that $G'/H'=G/H$. This proves (d2). In order to prove
(d3), let us assume first that $m=2$ and therefore $G \subset
S_1\times S_2$. Let us define
$$\phi':G\subset S_1\times S_2 \to S_2$$
as the restriction to $G$ of the projection map
    $S_1\times S_2\to S_2$.
Applying (d2) to $\pi'=\phi'$, we conclude that either
$\phi'(G)/\phi'(H) \cong G/H$ and one could put
$$\bar{H}:=\phi'(H) \subset \bar{G}:=\phi'(G) \subset
S_2$$ or
$$\bar{H}=\ker(\phi')\bigcap H \subset \bar{G}:=\ker{\phi'}\subset S_1
\times \{1\} \cong S_1.$$ In order to do the case of $m>2$ one has
only to write down
$$S_1\times \cdots S_m=S_1 \times (S_2 \times \cdots \times S_m)$$
and apply induction.
\item[(e)]
It is known (\cite{Margulis}, Remarks 1--2 on pp. 41-42; see also
\cite{Springer}, Ex. 16.3.9(1) on p. 277) that $\alpha(\K_1(E))$
is a normal subgroup in $\K_2(E)$ and the quotient
$\K_2(E)/\alpha(\K_1(E))$ is a commutative group of finite
exponent dividing $\deg(\alpha)$. Clearly, $G_{\alpha}$ coincides
with the kernel of the composition
$$G \stackrel{\gamma}{\longrightarrow}  \K_2(E)\twoheadrightarrow
\K_2(E)/\alpha(\K_1(E))$$ and, thanks to (i), must be an open
normal subgroup of finite index. Applying (ii)(a) to
$\pi=\alpha'$, we conclude that $H_{\alpha}$ is an open subgroup
of finite index in $G_{\alpha}$ and $G_{\alpha}/H_{\alpha}=G/H$.
Also, applying (ii)(a) to $G_{\alpha}\subset G
\stackrel{\gamma}{\longrightarrow} \K_2(E)$, we conclude that
$\gamma(G_{\alpha})$ is a compact subgroup of $\K_2(E)$,
$\gamma(H_{\alpha})$ is an open normal subgroup of finite index in
$\gamma(G_{\alpha})$ and
$\gamma(G_{\alpha})/\gamma(H_{\alpha})=G/H$.
 Since $\alpha$ is a
central isogeny, the kernel of $\alpha:\K_1(E) \to \K_2(E)$ is a
finite commutative group, whose exponent divides $\deg(\alpha)$.
This implies that the preimage $G_{-1,\alpha}$ of compact
$\gamma(G_{\alpha})$ in $\K_1(E)$ is also compact. Notice that
$$\alpha(G_{-1,\alpha})=\gamma(G_{\alpha}), \quad
\alpha(H_{-1,\alpha})=\gamma(H_{\alpha}).$$ Applying (ii)(a) to
$G_{-1,\alpha}\subset \K_1(E) \stackrel{\alpha}{\longrightarrow}
\K_2(E)$, we conclude that $H_{-1,\alpha}$ is an open normal
subgroup of finite index in $G_{-1,\alpha}$ and
$G_{-1,\alpha}/H_{-1,\alpha}=G/H$. In order to prove (e3) it
suffices to notice that the kernel of
$H_{-1,\alpha}\stackrel{\alpha}{\longrightarrow}\gamma(H_{\alpha})$
lies in the kernel of $\alpha:\K_1(E) \to \K_2(E)$ and therefore
also has exponent dividing $\deg(\alpha)$.
\end{enumerate}
\end{proof}
\begin{prop}
\label{proj} Suppose $V$ is a finite-dimensional vector space over
a finite field $k$ of characteristic $\ell$ and $G$ is a subgroup
of $\Aut(V)$ enjoying the following properties:
\begin{enumerate}
\item[(i)]
$G$ is perfect, i.e. $G=[G,G]$;
\item[(ii)]
 $G$ contains a normal abelian
subgroup $Z$ such that the quotient $\Gamma:=G/Z$ is a simple
non-abelian group.
\item[(iii)]
There exists a positive integer $d \ge \dim_k(V)$ such that every
nontrivial projective representation of $\Gamma$ in characteristic
$\ell$ has dimension $\ge d$.
\end{enumerate}
Then:
\begin{enumerate}
\item[(a)]
The $G$-module $V$ is absolutely simple, $\dim_k(V)=d$ and $Z$ is
the center of $G$. In particular, $Z$ consists of scalar matrices
and therefore is a cyclic group of order prime to $\ell$;
\item[(b)]
Every subgroup of $G$ (except $G$ itself) has index $\ge
\max(5,d+1)$;
\item[(c)]
For each finite field $k'$ of characteristic $\ell$ and each
positive integer $a<d$  every homomorphism $G \to \PGL_a(k')$ is
trivial.
\end{enumerate}
\end{prop}

\begin{proof}
{\bf Step 0}.  Each normal subgroup $H$ of $G$ either lies in $Z$
or coincides with $G$. Indeed, if $H$ contains $Z$ then the
simplicity of $G/Z$ implies that either $H=Z$ or $H=G$. Assume
that $H$ neither contains nor is contained in $Z$. Let us put
$D:=H\bigcap Z \ne Z$. Then $G_0=G/D$ is perfect, $Z_0=Z/D$ is a
nontrivial abelian normal subgroup in $G/D=G_0$ and
$G_0/H_0=G/Z=\Gamma$ is simple non-abelian and $H_0:=H/D$ is a
nontrivial normal subgroup of $G_0$ which meets $Z_0$ only at the
identity element. The simplicity of $G_0/Z_0$ implies that
$G_0=Z_0\times H_0$ which contradicts the perfectness of $G_0$.

We write $Z'$ for the Sylow-$\ell$-subgroup of $Z$. Clearly, $Z'$
is normal in $G$. It is also clear that if $Z'\ne\{1\}$ then every
semisimple faithful finite-dimensional representation of $Z'$ in
characteristic $\ell$ must have dimension $1$. On the other hand,
since $G$ is non-abelian, $\dim_k(V)>1$. This implies that the
$Z'$-module $V$ is semisimple if and only if $Z'=\{1\}$. Taking
into account that $Z'$ is normal in $G$ and applying Clifford's
theorem (\cite{CR}, \S 49, Th. 49.2), we conclude that if the
$G$-module $V$ is semisimple then $Z'=\{1\}$ and therefore the
order of $Z$ is not divisible by $\ell$.

 {\bf Step 1}. Assume that the $G$-module $V$ is
absolutely simple. Then $\#(Z)$ is prime to $\ell$. Let $k_1$ be
the finite field obtained by adjoining to $k$ all $\#(Z)$th roots
of unity and consider the $k_1$-vector space $V_1=V\otimes_k k_1$.
Clearly, $V_1$ carries a natural structure of absolutely simple
faithful $G$-module and
$$\dim_{k_1}(V_1)=\dim_k(V) \le d.$$
For each character $\chi: Z \to k_1^*$ we write $V^{\chi}$ for the
subspace
$$V^{\chi}:=\{v\in V_1\mid zv=\chi(z)v \quad \forall z\in Z \in
G\}.$$ Clearly,
$$V_1=\oplus_{\chi}V^{\chi},$$
 $G$ permutes all $V^{\chi}$'s and this action factors through
$G/Z=\Gamma$. It is also clear that the set of non-zero
$V^{\chi}$'s consists, at most, of $\dim_{k_1}(V_1)$ elements.
This implies that if the action of $G$ on all $V^{\chi}$'s is {\sl
non-trivial} then $G/Z=\Gamma$ contains a subgroup $S'\ne S$ with
index $r\le \dim_{k_1}(V_1) \le d$. This gives us a nontrivial
homomorphism $\Gamma \to \SS_r$ which must be an embedding, in
light of simplicity of $\Gamma$. This implies that $r\ge 5$ and
therefore $\SS_r$ is isomorphic to a subgroup of
$\PGL_{r-1}(\F_{\ell})$ and therefore $\Gamma$ is isomorphic to a
subgroup of $\PGL_{r-1}(\F_{\ell})$. Since $r\le d$, we get a
contradiction to property (iii). Hence, $G$ maps each $V^{\chi}$
into itself and therefore each $V^{\chi}$ is a $G$-invariant
subspace of $V_1$. The absolute simplicity of $V_1$ implies that
$V_1=V^{\chi}$ for some $\chi$. This implies that $Z\subset
k_1^*\I$; in particular, $Z$ is a central cyclic subgroup of $G$.
Since $G/Z$ is simple non-abelian, $Z$ coincides with the center
of $G$. Now the absolute simplicity of $V$ implies that $Z\subset
k^*\I\subset\Aut_k(V)$ and we get an embedding
$$\Gamma =G/Z\hookrightarrow \PGL(V).$$
This implies that $d\le \dim_k(V)$. Since $d \ge \dim_k(V)$, we
conclude that $d=\dim_k(V)$. This ends the proof of (a) in the
case of absolutely simple $G$-module $V$.

{\bf Step 2}. Assume that the $G$-module $V$ is simple (but not
necessarily absolutely simple). Let us put $\kappa=\End_G(V)$.
Clearly, $\kappa \supset k$ is a finite field of characteristic
$\ell$, $V$ carries a natural structure of absolutely simple
$\kappa[G]$-module and
$$\dim_{\kappa}(V)=\frac{\dim_k(V)}{[\kappa:k]}.$$
Applying the (special case of the) assertion (a) (proven in Step
1) to the absolutely simple $\kappa[G]$-module $V$, we conclude
that
$$d=\dim_{\kappa}(V)=\frac{\dim_k(V)}{[\kappa:k]} \le \dim_k(V)\le d.$$
We conclude that $\dim_k(V)=\dim_{\kappa}(V)$ and therefore
$[\kappa:k]=1$. This implies that $\kappa=k$, i.e., the $G$-module
$V$ is absolutely simple.  This ends the proof of (a) in the case
of  simple $G$-module $V$.

{\bf Step 3}. Assume that the $G$-module $V$ is semisimple (but
not necessarily  simple). Since $V$ is faithful, there is a simple
nontrivial $G$-submodule $W$ of $V$. We have
$$\dim_k(W) \le \dim_k(V) \le d.$$
Let us denote by $G_0 \ne\{1\}$ the image of $G$ in $\Aut_k(W)$.
Since $G$ is perfect, its homomorphic image $G_0$ is also perfect.
Since $G_0 \ne\{1\}$, the kernel of $G \twoheadrightarrow G_0$
lies in $Z$. We write $Z_0$ for the image of $Z$ in $G_0$.
Clearly, $Z_0$ is an abelian normal subgroup of $G_0$ and
$G_0/Z_0=G/Z=\Gamma$. Applying the (special case of the) assertion
(a) (proven in Step 2) to the faithful simple $G_0$-module $V$, we
conclude that
$$d=\dim_k(W)\le \dim_k(V) \le d.$$
This implies that $\dim_k(W)= \dim_k(V)$ and therefore $V=W$ is a
simple $G$-module. This ends the proof of (a) in the case of
semisimple $G$-module $V$.

 {\bf Step 4}.  End of the proof of (a). Assume that the $G$-module $V$
is {\sl not} semisimple. Let $V^{ss}$ be  its semisimplification.
Clearly, the $G$-module $V^{ss}$ is semisimple but {\sl not}
simple. Clearly, the kernel of the natural homomorphism $G \to
\Aut_k(V^{ss})$ consists of unipotent matrices and therefore is a
finite normal $\ell$-group. This implies that this kernel lies in
$Z$. Let $G_1$ be the image of $G$ in $\Aut_k(V^{ss})$ and $Z_1$
be the image of $Z$ in $G_1$. Clearly, $Z_1$ is an abelian normal
subgroup of $G_1$ and $G_0/Z_0=G/Z=\Gamma$. It is also clear that
the $G$-module $V^{ss}$ is faithful semisimple but not simple.
Applying the (special case of the) assertion (a) (proven in Step
3) to the faithful semisimple $G_1$-module $V^{ss}$, we conclude
that $V^{ss}$ is simple. We get a contradiction which proves that
$V$ is semisimple. This ends the proof of (a).

{\bf Step 5}. Proof of (c). Suppose $a<d$ is a positive integer,
$k'$ is a finite field of characteristic $\ell$ and $\phi: G \to
\PGL_a(k')$ is a nontrivial group homomorphism. By Step 0,
$\ker(\phi) \subset Z$. Let us put $G_2:=\phi(G)\subset
\PGL_a(k')$ and
$$Z_2:=\phi(Z) \subset G_2=\phi(G)\subset \PGL_a(k').$$
Clearly, $G_2$ is perfect, $Z_2$ is a central cyclic subgroup of
$G_2$ and $G_2/Z_2=\Gamma$. Let us denote by $G_3$ (resp. by
$Z_3$) the preimage of $G_2$ (resp. of $Z_2$) in $\GL_a(k')$ with
respect to the projectivization map $\GL_a(k')\to\PGL_a(k')$.
Clearly,
$$k_1^*\I\subset Z_3\subset G_3\subset \GL_a(k')$$
and
$$Z_3/k_1^*\I=Z_2 \subset G_2=G_3/k_1^*\I, \quad
G_3/Z_3=G_2/Z_2=\Gamma.$$

Since $k_1^*\I$ lies in the center of $Z_3$ (and of $G_3$) and the
quotient $Z_3/k_1^*\I=Z_2 $ is cyclic, $Z_3$ is abelian. If $G_3$
were perfect then we could apply the already proven assertion (a)
to the faithful $k'[G_3]$-module $k{'}^a$ and conclude that
$a=\dim_{k'}k{'}^a=d$. This would lead us to contradiction, since
$a<d$ and we conclude that $\phi$ is trivial and we are done.
However, there is no reason to believe that $G_3$ is perfect. So,
let us choose a minimal  subgroup $G'$ of $G_3$ which maps {\sl
onto} $G_3/Z_3=\Gamma$. Such a a choice is possible in light of
finiteness of $G_3$. Let us denote by $Z'$ the intersection of
$G'$ and $Z_3$. Clearly, $Z'$ is an abelian normal subgroup in
$G'$ and $G'/Z'=G_3/Z_3=\Gamma$. The minimality of $G'$ and
perfectness of $\Gamma$ imply that $G'$ is perfect.
 Now, applying the assertion (a) to
  the faithful $k'[G']$-module $k{'}^a$, we conclude that
$a=\dim_{k'}k{'}^a=d$. Since $a<d$ our assumption that $\phi$ is
non-trivial was wrong. This proves (c).

{\bf Step 6}. We still have to prove (b). Assume that $G$ contains
a subgroup of index $<5$. This gives us a non-trivial homomorphism
of $G$ into the solvable group $\SS_4$. Since its kernel must lie
in the abelian group $Z$, we conclude that $G$ is solvable which
is not the case. So, if $d<5$ then we are done.

 Assume that $d\ge 5$ and $G$ contains a subgroup of index $r\le d$. This gives
us a non-trivial homomorphism of $G$ into $\SS_r\subset\SS_d$.
Since $\SS_d$ is isomorphic to a subgroup $\PGL_{d-1}(\F_{\ell})$,
we conclude that there exists a nontrivial homomorphism $G \to
\PGL_{d-1}(\F_{\ell})$. But this contradicts (c).
\end{proof}

\begin{rem}
\label{repA}
 Let $d \ge 8$ be an even integer.
 \begin{enumerate}
 \item[(a)]
 In characteristic $2$ all nontrivial projective
representations of $\A_{d+1}$ and of $\A_{d+2}$ have dimension
$\ge d$. Indeed, it is well-known that the groups $\A_{d+1}$ and
$\A_{d+2}$ are perfect and their Schur multipliers coincide and
equal to $2$. This implies that all irreducible projective linear
representation of $\A_{d+1}$ and of $\A_{d+2}$ in characteristic
$2$ are, in fact, linear representations. By a theorem of Wagner
\cite{Wagner}, all nontrivial linear representation of $\A_{d+1}$
and of $\A_{d+2}$ in characteristic $2$ have dimension $\ge d$.
This implies that all irreducible projective  representation of
$\A_{d+1}$ and of $\A_{d+2}$ in characteristic $2$ have dimension
$\ge d$. Taking into account that $\A_{d+1}$ and $\A_{d+2}$ are
simple non-abelian groups, we conclude that all nontrivial
projective  representation of $\A_{d+1}$ and of $\A_{d+2}$ in
characteristic $2$ have dimension $\ge d$.
\item[(b)]
Each nontrivial projective representation of $\A_{d+2}$ in
characteristic zero has dimension $\ne d$. Indeed, $d+2$ is an
even integer $\ge 10$ and the desired assertion about
representations in characteristic zero was proven in
\cite{ZarhinMRL2}.
\item[(c)]
There does {\sl not} exist a faithful symplectic $d$-dimensional
representation of $\A_{d+1}$ in characteristic $0$. Indeed, there
exists exactly one (up to an isomorphism) faithful $d$-dimensional
representation of $\A_{d+1}$ in characteristic $0$ (Th. 2.5.15 on
p. 71 of \cite{JK}) and this representation is absolutely
irreducible and orthogonal; hence it could not be symplectic.
\item[(d)]
Suppose a short exact sequence
$$1 \to \Z/2\Z \hookrightarrow \A_{d+1}' \twoheadrightarrow
\A_{d+1}\to 1$$ defines a non-splitting central extension of
$\A_{d+1}$ (i.e., $\A_{d+1}'$ is the universal central extension
of $\A_{d+1}$). If $d \ge 10$ then every faithful representation
of $\A_{d+1}'$ in characteristic $0$ has dimension $\ne d$.
Indeed, for $10 \le d \le 12$ this assertion follows from the
character tables in \cite{HH}. So, further we assume that $d \ge
14$.

We start with an elementary discussion of the dyadic expansion
$d+1=2^{w_1}+ \cdots + 2^{w_s}$ of $d+1$. Here $w_i$'s are
distinct nonnegative integers with $w_1< \cdots  < w_s$ and $s$ is
the exact number of terms (non-zero digits) in the dyadic
expansion of $n$. Since $d+1$ is odd, $w_0=0$ and  each $w_i \ge
i-1$. This implies that $d+1 \ge 2^s-1$ and therefore
$$s \le \log_2(d+2).$$

By a theorem of Wagner (Th. 1.3(ii) on pp. 583--584 of
\cite{Wagner2}), each proper projective representation of
$\A_{d+1}$ (i.e., a nontrivial linear representation of
$\A_{d+1}'$) in characteristic $\ne 2$ has dimension divisible by
$N:=2^{\lfloor\frac{(d+1)-s-1}{2}\rfloor}=2^{\lfloor\frac{(d-s}{2}\rfloor}$.
So, in order to prove (b), it suffices to check that $d$ is {\sl
not} divisible by $N$ for all even $d\geq 14$.

 If $d\geq 14$ then $2^{d-1}> (d+2)d^2$. Then
$2^{d-\log_2(d+2)-1}>d^2$ and therefore  $2^{d-s-1}>d^2$. Taking
square roots at both sides, we get $2^{\frac{d-s-1}{2}}>d$. Then
we see easily that $N=2^{\lfloor\frac{d-s}{2}\rfloor}>d$. This
finishes the proof.

\item[(d1)]
If $d=8$ then all faithful  absolutely irreducible representations
of $\A_{d+1}'=\A_{9}'$ in characteristic $0$ have dimension $\ge
8$; among them all the $8$-dimensional representations are
orthogonal (\cite{Atlas}, p. 37). As above, this implies that
there does {\sl not} exist a faithful symplectic $d$-dimensional
representation of $\A_{d+1}'$ in characteristic $0$.
\end{enumerate}
\end{rem}

\begin{rem}
\label{repM}
Let us put $d=10$. Recall \cite{Atlas} that the Schur
multiplier of  the Mathieu group $\M_{11}$ is $1$ and therefore
all projective representations of $\M_{11}$ are, in fact, linear.
It is known \cite{AtlasB} that all faithful irreducible
representation of $\M_{11}$ in characteristic $2$ have dimension
$\ge 10$. Since $\M_{11}$ is perfect, all nontrivial
representation of $\M_{11}$ in characteristic $2$ have dimension
$\ge 10$. It is also known \cite{Atlas} that in characteristic $0$
all faithful irreducible representation of $\M_{11}$  have
dimension $\ge 10$ and none of $10$-dimensional absolutely
irreducible representations of $\M_{11}$ is symplectic. This
implies that in characteristic $0$ none of faithful
$10$-dimensional representations of $\M_{11}$ is symplectic.
\end{rem}

\begin{rem}
\label{repL} Let $q$ be an {\sl odd} power prime, $m \ge 3$ a
positive integer, $B=\P^{m-1}(\F_q)$ the $(m-1)$-dimensional
projective space over $\F_q$. Clearly, the cardinality of $B$ is
$\frac{q^m-1}{q-1}$ which is odd (resp. even) if $m$ is odd (resp.
even). The projective special linear group $\L_m(q):=\PSL(m,\F_q)$
acts naturally, faithfully  and doubly transitively on
$B=\P^{m-1}(\F_q)$. It is well known that this action gives rise
to the {\sl deleted permutation representation}: a certain
faithful absolutely irreducible $\Q[\L_m(q)]$-module $(\Q^B)^0$ of
$\Q$-dimension $\frac{q^m-1}{q-1}-1$ (see for instance
\cite{ZarhinTexel}). Clearly, the representation of $\L_m(q)$ in
$(\Q^B)^0$ is orthogonal (since it is defined over $\Q$) and
therefore is {\sl not} symplectic, because it is absolutely
irreducible.  It is known (\cite{TZ}, Th. 1.1) that if $(q,n) \ne
(3,4)$ then in characteristic $0$ all nontrivial irreducible
projective representations of $\L_m(q)$ have dimension $\ge
\frac{q^m-1}{q-1}-1$ and there is exactly one (up to an
isomorphism) a nontrivial irreducible projective representation of
$\L_m(q)$ of dimension $\frac{q^m-1}{q-1}-1$. This implies easily
that if a short exact sequence
$$1 \to \Z/2\Z \hookrightarrow \L_m(q)' \twoheadrightarrow
\L_m(q)\to 1$$ defines a central extension of $\L_m(q)$ then in
characteristic $0$ there  does not exist a faithful symplectic
absolutely irreducible representation of $\L_m(q)'$ of dimension
$\le \frac{q^m-1}{q-1}-1$.

 Guralnick  proved that if $(q,n) \ne (3,4)$ then the
dimension of each nontrivial projective irreducible representation
of $\L_m(q)$ in characteristic $2$ is greater than or equal to
$$2[(\frac{q^{m}-1}{q-1}-1)/2]$$
(see Th. 1.1 and Table III in \cite{GurTiep}). Since $\L_m(q)$,
this implies easily that in in characteristic $2$ the dimension of
every nontrivial projective  representation of $\L_m(q)$ in
characteristic $2$ is greater than or equal to
$$2[(\frac{q^{m}-1}{q-1}-1)/2].$$
\end{rem}

\section{Very simple representations}
The following notion was introduced by the author in
\cite{ZarhinTexel}.

\begin{defn}
Let $V$ be a vector space over a field $k$, let $G$ be a group and
$\rho: G \to \Aut_{k}(V)$ a linear representation of $G$ in $V$.
We say that the $G$-module $V$ is {\sl very simple} if it enjoys
the following property:

If $R \subset \End_{k}(V)$ is an $k$-subalgebra containing the
identity operator $\I$ such that

 $$\rho(\sigma) R \rho(\sigma)^{-1} \subset R \quad \forall \sigma \in G$$
 then either $R=k\cdot \I$ or $R=\End_{k}(V)$.

 Here is (obviously) an equivalent definition: if $R \subset \End_{k}(V)$ is an $k$-subalgebra containing the
identity operator $\I$ and stable under the conjugations by all
$\rho(\sigma)$ then either $\dim_k(R)=1$ or
$\dim_k(R)=(\dim_k(V))^2$.
\end{defn}

\begin{rems}
\label{image}
\begin{enumerate}
\item[(i)]
Clearly, the $G$-module $V$ is very simple if and only if the
corresponding $\rho(G)$-module $V$ is very simple.
\item[(ii)]
Clearly, if $V$ is very simple then the corresponding algebra
homomorphism
    $$k[G] \to \End_{k}(V)$$
is surjective. Here $k[G]$ stands for the group algebra of $G$. In
particular, a very simple module is absolutely simple.
\item[(iii)]
If $G'$ is a subgroup of $G$ and the $G'$-module $V$ is very
simple then the $G$-module $V$ is also very simple.
\item[(iv)]
Let $G'$ be a normal subgroup of $G$. If $V$ is a very simple
$G$-module then either $\rho(G') \subset \Aut_{k}(V)$ consists of
scalars (i.e., lies in $k\cdot\I$) or the $G'$-module $V$ is
absolutely simple. Indeed, let $R'\subset \End_{k}(V)$ be the
image of the natural homomorphism $k[G'] \to \End_{k}(V)$.
Clearly, $R'$ is stable under the conjugation by elements of $G$.
Hence either $R'$ consists of scalars and therefore
$\rho(G')\subset R'$ consists of scalars or $R'= \End_{k}(V)$ and
therefore the $G'$-module $V$ is absolutely simple.
\item[(v)]
Suppose $F$ is a discrete valuation field with  valuation ring
$O_F$, maximal ideal $m_F$ and residue field $k=O_F/m_F$. Suppose
$V_F$ a finite-dimensional $F$-vector space,
$$\rho_F: G \to \Aut_F(V_F)$$
a $F$-linear representation of $G$. Suppose $T$ is a $G$-stable
$O_F$-lattice in $V_F$ and the corresponding $k[G]$-module $T/m_F
T$ is isomorphic to $V$. Assume that the $G$-module $V$ is very
simple. Then:
\begin{enumerate}
\item[(a)]
The $G$-module $V_F$ is also very simple. In other words, a
lifting of a very simple module is also very simple. Indeed, let
$R_F \subset \End_F(V_F)$ be an $F$-subalgebra containing the
identity operator and stable under the conjugation by elements of
$G$. Let us put
$$R_O=R \bigcap \End_{O_F}(T) \subset \End_{O_F}(T).$$
Clearly, $R_O$ is a free $O_F$-module, whose rank coincides with
$\dim_F(R_F)$. It is also clear that $R_O$ is a pure
$O_F$-submodule of $\End_{O_F}(T)$. This implies that
$$R_{k}=R_O/m_FR_O=R_O\otimes_{O_F}k\subset\End_{O_F}(T)\otimes_{O_F}k=
\End_{k}(T/m_F T)=\End_{k}(V)$$ is an $k$-subalgebra of
$\End_{k}(V)$ of dimension $\dim_F(R_F)$, contains the identity
operator and is stable under the conjugation by elements of $G$.
Now the very simplicity of $V$ implies that either
$\dim_{k}(R_{k})=1$ or $\dim_{k}(R_{k})=\dim_{k}(V)^2$. Since
$\dim_{k}(V)=\dim_F(V_F)$, we conclude that either $\dim_F(R_F)=1$
or $\dim_F(R_F)=\dim_F(W)^2$. Clearly, in the former case $R_F$
consists of scalars and in the latter one $R_F=\End_F(V_F)$.
\item[(b)]
Suppose that $\fchar(F)=0$ and $\rho_F$ is an embedding. Further
we identify $G$ with its image $\rho_F(G)\subset \Aut(V)$.

As usual, we write $\GL_{V_F}$ for the $F$-algebraic group of
automorphisms of $V$. In particular, $\GL_{V_F}(F)=\Aut(V)$. We
write $\SL_{V_F}$ for the $F$-algebraic group of automorphisms of
$V$ with determinant $1$. In particular, $\SL_{V_F}(F)$ coincides
with $\SL(V)$.

Let $G_{\alg}$ be the {\sl algebraic envelope of G}. i.e., the
smallest algebraic subgroup of $\GL_{V_F}$, whose group of
$F$-points contains $G$. Since $G\subset G_{\alg}(F)$, the
$G_{\alg}(F)$-module $V_F$ is also very simple. Clearly, if
$G\subset \SL(V_F)$ then $G_{\alg}\subset \SL_{V_F}$.

Let $G_{\alg}^0$ be the identity component of $G_{\alg}$. Clearly,
$G_{\alg}^0(F)$ is a normal subgroup of finite index in
$G_{\alg}(F)$. By (iv), either $G_{\alg}^0(F)$ consists of scalars
or the $G_{\alg}^0(F)$-module $V_F$ is absolutely simple.

Assume, in addition, that if $G\subset \Aut_F(V_F)$ is infinite
and lies in $\SL(V_F)$ then $G_{\alg}(F)$ is also infinite and
lies in $\SL(V_F)$. Since $G_{\alg}^0(F)$ has finite index in
$G_{\alg}(F)$, we conclude that $G_{\alg}^0(F)$ is also infinite
and lies in $\SL(V_F)$. This implies that $G_{\alg}^0(F)$ could
not consist of scalars and therefore the $G_{\alg}^0(F)$-module
$V_F$ is absolutely simple. This implies easily that the
$F$-algebraic group $G_{\alg}^0$ is semisimple.
\item[(c)]
Suppose $\ell$ is a prime, $F=\Q_{\ell}, O_F=\Z_{\ell},
k=\F_{\ell}$. We write $V_{\ell}$ for $V_F=V_{\Q_{\ell}}$.
 Assume that $G \subset \Aut(V_{\ell})$ is a
compact subgroup. Then $G$ is a compact $\ell$-adic Lie subgroup
of $\Aut(V_{\ell})$ (\cite{SerreLie}). We write $\g$ for its Lie
algebra; it is a $\Q_{\ell}$-Lie subalgebra of $\End(V_{\ell})$.
The absolute simplicity of the $G$-module $V_{\ell}$ implies that
the $\g$-module $V_{\ell}$ is semisimple (Prop. 1 of
\cite{SerreLocal}) and therefore $\g$ is reductive, i.e.,
$$\g=\c \times \s$$
where $\c$ is the center of $\g$ and $\s=[\g,\g]$ is a semisimple
$\Q_{\ell}$-Lie algebra. This implies also that that the
$F$-algebraic group $G_{\alg}^0$ is reductive (Prop. 2 of
\cite{SerreLocal}) and its Lie algebra $\Lie(G_{\alg}^0)$
coincides with $\c_{\alg}\times s$ where $\c_{\alg}$ is the
algebraic envelope of $\c$  and coincides with the center of
$\Lie(G_{\alg})^0)$.

 We keep the assumption
that $G\subset \Aut_{\Q_{\ell}}(V_{\ell})$ is infinite and lies in
$\SL(V_{\ell})$. It follows that that $\Lie(G_{\alg}^0)$ lies in
the Lie algebra $\sll(V_{\ell})$ of linear operators in $V_{\ell}$
with zero trace. Since the $G_{\alg}^0(F)$-module $V_{\ell}$ is
absolutely simple (by (b)), the natural representation of the
connected $\Q_{\ell}$-algebraic group $G_{\alg}^0$ in $V_{\ell}$
is also absolutely irreducible. This implies, in turn, that the
natural representation of the Lie algebra $\Lie(G_{\alg}^0)$ in
$V_{\ell}$ is also absolutely irreducible and therefore its center
$\c_{\alg}$ is either zero or consists of scalars. Since
$$\c_{\alg} \subset\c_{\alg}\times s=\Lie(G_{\alg}^0)\subset \sll(V_{\ell}),$$
we conclude that there are no non-zero scalars in $\c_{\alg}$ and
therefore $\c_{\alg}=\{0\}$. This implies that $\c=\{0\}$ and
therefore
$$\Lie(G_{\alg}^0)=\s=\g.$$
In particular, $\g$ is a semisimple algebraic $\Q_{\ell}$-Lie
algebra. This implies that $G$ is open in $G_{\alg}(\Q_{\ell})$ in
$\ell$-adic topology and meets all the components of $G_{\alg}$
(\cite{SerreLocal}, Prop. 2 and its Corollary).
\end{enumerate}
\end{enumerate}
\end{rems}

The following assertion is a special case of Th. 4.3 of
\cite{ZarhinCrelle}.

\begin{thm}
\label{veryfactor} Suppose a finite field field $k$, a positive
integer $N$ and a group $H$ enjoy the following properties:

\begin{itemize}
 \item
 $H$ is perfect, i.e., $H=[H,H]$;
\item
Each homomorphism from $H$ to $\SS_N$ is trivial;
\item
Let $N=ab$ be a factorization of $N$ into a product of two
positive integers $a$ and $b$. Then either each homomorphism from
$H$ to $\PGL_a(k)$ is trivial or each homomorphism from $H$ to
$\PGL_b(F)$ is trivial.
\end{itemize}
Then each absolutely simple $H$-module of $k$-dimension $N$ is
very simple. In other words, in dimension $N$ the properties of
absolute simplicity and very simplicity over $F$ are equivalent.
\end{thm}

\begin{cor}
\label{verymin} Suppose $V$ is a finite-dimensional vector space
over a finite field $k$ of characteristic $\ell$ and $G$ is a
subgroup of $\Aut(V)$ enjoying the following properties:
\begin{enumerate}
\item[(i)]
$G$ is perfect, i.e. $G=[G,G]$;
\item[(ii)]
 $G$ contains a normal abelian
subgroup $Z$ such that the quotient $\Gamma:=G/Z$ is a simple
non-abelian group.
\item[(iii)]
There exists a positive integer $d \ge \dim_k(V)$ such that every
nontrivial projective representation of $\Gamma$ in characteristic
$\ell$ has dimension $\ge d$.
\end{enumerate}
Then $Z$ is a cyclic central subgroup of $G$ and the $G$-module
$V$ is very simple.
\end{cor}

\begin{proof}
Replacing $G$ by its minimal subgroup which maps onto $\Gamma$ we
may assume, in light of Remark \ref{image}(iii) that $G$ is
perfect. Now the very simplicity follows readily from Prop.
\ref{proj} combined with Th. \ref{veryfactor}.
\end{proof}

\section{$\ell$-adic representations}
\begin{thm}
\label{Lie} Suppose $\ell$ is a prime, $V_{\ell}$ is a
$\Q_{\ell}$-vector space of finite dimension $d>1$, $G \subset
\SL(V_{\ell})\subset \Aut(V_{\ell})$ a compact linear group. We
write $G_{alg}$ for the algebraic envelope of $G$ and $\g$ for the
Lie algebra of $G$.

 Suppose $T\subset V_{\ell}$ is a $G$-stable
$\Z_{\ell}$-lattice in $V_{\ell}$. Let us put $V(\ell):=T/\ell T$.
Clearly, $V(\ell)$ is a $d$-dimensional vector space and carries a
natural structure of $G$-module. Let us denote by
$\tilde{G}_{\ell}$ the image of the natural homomorphism
$$G \to \Aut_{\Z_{\ell}}(T) \twoheadrightarrow \Aut(T/\ell T)=\Aut_{\F_{\ell}}(V(\ell)).$$
Suppose  $\tilde{G}_{\ell}$ enjoys the following properties:
\begin{enumerate}
\item[(a)]
$\tilde{G}_{\ell}$ is a simple non-abelian group;
\item[(b)]
Every faithful
 projective representation of $\tilde{G}_{\ell}$ in
characteristic $\ell$ has dimension $\ge d$;
\item[(c)]
Either $\ell=2$ or $\tilde{G}_{\ell}$ is a finite group of Lie
type in odd characteristic or one of $26$ known sporadic groups.
(In other words, modulo the classification, either $\ell=2$ or
$\tilde{G}_{\ell}$ is not a group of Lie type in characteristic
$2$.)
\item[(d)]
Either $G$ is infinite or one the two following conditions holds:
\begin{enumerate}
\item[(i)]
$\ell>2$ and there does not exist a lifting of the
$\F_{\ell}[\tilde{G}_{\ell}]$-module $V(\ell)$ to an absolutely
simple $\Q_{\ell}[\tilde{G}_{\ell}]$-module $\Q_{\ell}^d$. (E.g.,
every homomorphism of $\tilde{G}_{\ell}\to \GL_d(\Q_{\ell})$ is
trivial.)
\item[(ii)]
$\ell=2$ and for each central extension (short exact sequence)
$$1 \to \Z/2\Z \hookrightarrow \tilde{G}' \twoheadrightarrow
\tilde{G}_{2}\to 1$$ there does not exist such a lifting of the
$\F_{2}[\tilde{G}']$-module $V(2)$ to an absolutely simple
$\Q_{2}[\tilde{G}']$-module $\Q_{2}^d$ that the distinguished
central subgroup $\Z/2\Z$ acts on $\Q_{2}^d$ via multiplications
by $\pm 1$. (E.g., every homomorphism  $\tilde{G}_{2}\to
\PGL_d(\Q_{2})$ is trivial.)
\item[(iii)]
$\ell=2$, there exists a non-degenerate $G$-invariant alternating
bilinear form on $V_{\ell}$ and for each central extension
$$1 \to \Z/2\Z \hookrightarrow \tilde{G}' \twoheadrightarrow
\tilde{G}_{2}\to 1$$ there does not exist such a lifting of the
$\F_{2}[\tilde{G}']$-module $V(2)$ to an absolutely simple
$\Q_{2}[\tilde{G}']$-module $\Q_{2}^d$ that the distinguished
central subgroup $\Z/2\Z$ acts on $\Q_{2}^d$ via multiplications
by $\pm 1$ and there exists a non-degenerate
$\tilde{G}'$-invariant alternating bilinear form on $\Q_{2}^d$.
\end{enumerate}

\end{enumerate}

Then $G$ is an open subgroup of $G_{alg}(\Q_{\ell})$ and the
$\Q_{\ell}$-Lie algebra $\g$ is  absolutely simple . If $E$ is a
finite Galois extension of $\Q_{\ell}$ such that
 the $E$-simple Lie algebra $g_E:=g\otimes_{\Q_{\ell}}E$ splits
 then the faithful $g_E$-module $V_E=V\otimes_{\Q_{\ell}}E$
 is a fundamental representation of minimum dimension.
\end{thm}

\begin{proof}

{\bf Step 0}. It follows from Corollary \ref{verymin} that the
$\tilde{G}_{\ell}$-module $V(\ell)$ is very simple. In turn, it
follows from Remark \ref{image}(v) that the $G$-module $V_{\ell}$
is very simple. Let us denote by $H$ the kernel of $G
\twoheadrightarrow \tilde{G}_{\ell}$. Clearly, $H$ is a closed
normal subgroup $G$ and, by Remark \ref{image}(iv), either $H$
consists of scalars or the $H$-module $V_{\ell}$ is absolutely
simple. If $G$ is finite then $H$ is also finite. Notice also that
$H$ is a pro-$\ell$-group. We have
$$G/H=\tilde{G}_{\ell}.$$
The simplicity of $\tilde{G}_{\ell}$ implies that $H$ is the
largest closed {\sl normal} pro-$\ell$-subgroup in $G$.

 {\bf Step 1}. $G$
is {\sl infinite}. Indeed, let us assume that $G$ is finite. Then
$H$ is a finite group consisting of automorphisms of $T$ congruent
to $1$ modulo $\ell$.

If $\ell=2$ then, thanks to Minkowski-Serre Lemma \cite{SZ}, all
nontrivial elements of $H$ have order $2$ and therefore $H$ is a
finite abelian group. Since every absolutely irreducible
representation of a finite abelian group must have dimension $1$,
we conclude that $H$ consists of scalars and therefore either
$H=\{1\}$ or $H=\{\pm 1\}$. In both cases if we denote by
$\tilde{G}'$ the subgroup of $\Aut_{\Z_2}(T)$ generated by $G$ and $\{\pm
1\}$ then the reduction map modulo $2$ gives us a central
extension
$$1 \to  \{\pm 1\}\hookrightarrow \tilde{G}' \twoheadrightarrow
\tilde{G}_{2}\to 1.$$ Clearly, the $\tilde{G}'$-module
$V_{2}=\Q_2^d$ is very simple and is a lifting of the
$\tilde{G}_{2}$-module $V_{2}$. It is also clear that the
$\tilde{G}'$-module $V_{2}=\Q_2^d$ is symplectic if the $G$-module
$V_{2}$ is symplectic.
 This contradicts (d)(ii) and (d)(iii) respectively and
therefore $G$ must be infinite.

If $\ell>2$ then, thanks to Minkowski-Serre Lemma \cite{SZ},
$H=\{1\}$ and the reduction map gives us an isomorphism $G \cong
\tilde{G}_{\ell}$. This implies easily that the
$G=\tilde{G}_{\ell}$-module $V_{\ell}\cong \Q_{\ell}^d$ is a
lifting of the $\tilde{G}_{\ell}$-module $V(\ell)$. This
contradicts (d)(i) and therefore $G$ must be infinite.

 {\bf Step 2}. Now we know that $G$ is infinite and lies
in $\SL(V_{\ell})$. It follows from Remark \ref{image}(v)(c) that
the identity component $G_{\alg}^0$ of $G_{\alg}$ is a semisimple
$\Q_{\ell}$-algebraic group, its Lie algebra coincides with $\g$
and $G$ is an open subgroup of $G_{\alg}(\Q_{\ell}$. In addition,
$\g$ is a semisimple $\Q_{\ell}$-Lie algebra coinciding with the
Lie algebra of $G_{\alg}$ and the natural representation of $\g$
in $V_{\ell}$ is absolutely irreducible.

{\bf Step 3}. There exists a finite Galois extension $E$ of
$\Q_{\ell}$ such that the semisimple $E$-Lie algebra
$$\g_E:=\g\otimes_{\Q_{\ell}}E$$
is split; in particular, $\g_E$ splits into a direct sum
$$\g_E=\oplus_{i \in I}\g_i$$
of absolutely simple split $E$-Lie algebras $\g_i$. Here $I$ is
the set of minimal non-zero ideals $\g_i$ in $\g_E$. It is
well-known that
$$V_E:=V_{\ell}\otimes_{\Q_{\ell}}E$$
becomes a faithful absolutely simple $\g_E$-module and splits into
a tensor product
$$V_E=\otimes_{i \in I} W_i$$
of faithful absolutely simple $\g_i$-modules $W_i$. Since each
$g_i$ is simple and $W_i$ is faithful,
$$\dim_{E}(W_i) \ge 2 \quad \forall i \in I.$$
This implies that the cardinality
$$r:=\#(I) \le \log_2(\dim_E(V_E))=\log_2(d)<d.$$

Let us consider the adjoint representation
$$\Ad: G \to \Aut(\g) \subset \Aut(\g_E).$$
Since the $\g$-module $V_{\ell}$ is absolutely simple, $\ker(\Ad)$
coincides with the finite subgroups of scalars in $G$. (The
finiteness folows from the inclusion $G \subset \SL(V_{\ell})$. It
follows easily that $\ker(\Ad)$ coincides with the center $\ZZ(G)$
of $G$.

 Clearly, $\Ad$ permutes elements of $I$ and therefore gives rise to the continuous
homomorphism (composition)
$$\pi_1:G \to \Aut(\g_E) \to \Perm(I) \cong \SS_r.$$
Clearly, one could embed $\SS_r$ into $\PGL_r(\bar{\F}_{\ell})$.
Since $r<d$, it follows from Proposition \ref{cont}(ii)(c) that
$\tilde{G}_{\ell}=G/H=G_1/H_1$ where
$$G_1=\ker(\pi_1)\subset \Aut_{\Z_{\ell}}(T), \quad  H_1=\ker(\pi_1)\bigcap H.$$
Clearly, $G_1$ is an open subgroup of finite index in $G$ and
therefore its Lie algebra coincides with $\g$. It is also clear
that $H_1$ is a pro-$\ell$-group. By definition of $G_1$ the image
of
$$G_1 \subset G \to \Aut(\g) \subset \Aut(\g_E)=\Aut(\oplus_{i \in I}\g_i)$$
lies in $\prod _{i \in I}\Aut(\g_i)$.
Let us put
$$\G=G_{\alg}^0.$$
We write $\G_E$ for the  semisimple split $E$-algebraic subgroup
of $\GL_{V_E}$ obtained from $\G$ by extension of scalars.
Clearly, the $E$-Lie algebra of $\G_E$ coincides with $\g_E$.

We write $\G_i$ for the {\sl simply connected} absolutely simple
split $E$-algebraic subgroup, whose Lie algebra coincides with
$\g_i$ (\cite{Springer}).

 We write
$\G_E^{\Ad} \subset \GL_{\g_E}$ for the adjoint group of $G_E$. If
$\G_i^{\Ad} \subset \GL_{\g_i}$ is the adjoint  group of $\G_i$
then
$$\G_E^{Ad}=\prod _{i \in I}\G_i^{\Ad}.$$
It is well-known that for each $i\in I$ the group $\G_i^{\Ad}(E)$
is a closed (in $\ell$-adic topology) normal subgroup in
$\Aut(\g_i)$ of finite index $1,2$ or $6$; in the latter case the
quotient $\Aut(\g_j)/\G_j^{\Ad}(E)$ is isomorphic to $\SS_3$.
(Recall that $\g_j$ is split.) Let us consider the composition
$$\pi_2:G_1 \to \prod _{i \in I}\Aut(\g_i) \twoheadrightarrow
\prod _{i \in I}\Aut(\g_i)/\prod _{i \in I}\G_i^{\Ad}(E)= \prod_{i
\in I} \Aut(\g_i)/\G_i^{\Ad}(E).$$

It follows from Proposition \ref{cont}(ii)(b) that if we put
$$G_2:=\ker(\pi_2)\subset G_1 \subset \Aut_{\Z_{\ell}}(T), \quad H_2=H_1\bigcap G_2$$
then $G_2$ is compact, $H_2$ is an open normal pro-$\ell$-subgroup
of finite index in $G_2$ and $G_2/H_2=\tilde{G}_{\ell}$.
 By definition of $G_2$ the image of
$$G_2 \subset G_1 \to \prod_{i \in I}\Aut(\g_i)\twoheadrightarrow \Aut(g_j)$$
lies in $\G_j^{\Ad}(E)$ for all $j \in I$.

{\bf Step 4}. Let us consider the  canonical central isogeny of
semisimple $E$-algebraic groups
$$\alpha:=\prod _{i \in I}\Ad_{i,E}:\prod _{i \in I}\G_i \to \prod _{i \in I}\G_i^{\Ad}=\G_E^{\Ad}.$$
We also have the continuous group homomorphism
$$\pi_2: G_2 \to \prod _{i \in I}\G_i^{\Ad}(E)=\G_E^{\Ad}(E),$$
whose kernel is a finite commutative group (consisting of
scalars). Applying Prop. \ref{cont}(ii)(e2), we conclude that
there exists is a compact subgroup $G_3 \subset \prod _{i \in
I}\G_i^{\Ad}(E)=\G_E(E)$ and an open normal subgroup $H_3$ of
finite index in $G_3$ such that $G_3/H_3 \cong \tilde{G}_{\ell}$.
By Prop. \ref{cont}(ii)(e3), every prime divisor of the
(super)order of $H_3$ is either $\ell$ or a divisor of one of
$\deg(\Ad_{i,E})$.
 Applying Prop. \ref{cont}(ii)(d3), we conclude that there exist
$j\in I$, a compact subgroup $G_4 \subset \G_j(E)$ and an open
normal subgroup $H_4$ of finite index in $G_4$ such that $G_4/H_4
\cong \tilde{G}_{\ell}$. In addition, every prime divisor of the
(super)order of $H_4$ is either $\ell$ or a divisor of one of
$\deg(\Ad_{i,E})$.

 {\bf Step 5}.
 Let $W$ be a finite-dimensional $E$-vector space which carries a
 structure of faithful absolutely simple $g_j$-module. I claim
 that
 $$d':=\dim_E(W) \ge d.$$
 In order to prove this inequality first, notice that there exists
 a $E$-rational representation
 $$\rho_W:\G_i \to \GL_W,$$
 whose kernel is a finite central subgroup of $\G_i$. Let us
 consider the  continuous homomorphism $\pi_4$ defined as
 $$ G_4 \to \G_j(E)\stackrel{\rho_W}{\longrightarrow}\Aut(W).$$
 Clearly, $\ker(\pi_4)$ is a finite commutative group. Applying
  Prop. \ref{cont}(ii)(d2), we conclude that $G_5=\rho_W(G_4)$ is
  a compact subgroup of $\Aut(W)$ containing an open normal
  subgroup $H_5:=\rho_W(H_4)$ and $\tilde{G}_{\ell}=G_5/H_5$.

 Since $G_5$ is compact, there is a $G_5$-stable $O_E$-lattice $T_O$
 in $W$. Notice that
 Here $O_E$ stands for the ring of
 integers in $E$. Let $\m_E$ be the maximal ideal in $O_E$ and
 $k=O_E/m_E$ be the corresponding residue field. Let us denote by
 $\pi_5$
 the restriction
 $$G_5 \subset \Aut_{O_E}(T) \twoheadrightarrow
 \Aut_{k}(T/\m_E)T)\cong \GL_{d'}(k).$$
 of the residue map to $G_5$.
 Clearly, $\ker(\pi_5)$ lies in the kernel of the reduction map
$\Aut_{O_E}(T) \twoheadrightarrow \Aut_{k}(T/\m_E)T)$ and
therefore is a pro-$\ell$-group. Hence $\tilde{G}_{\ell}$ could
not be a homomorphic image of $\ker(\pi_5)$. Let us put
$$G_6=\pi_5(G_5) \subset \GL_{d'}(k), \quad H_6=\pi_5(H_5) \subset
G'_6.$$ Applying Proposition \ref{cont}(ii)(d2) to $G_5$ and
$H_5$, we conclude that
$$\tilde{G}_{\ell}=G_5/H_5=G_6/H_6.$$
Let us consider the projectivization map

$$\pi_6:G_6\subset \GL_{d'}(k) \twoheadrightarrow
\PGL_{d'}(k).$$ Clearly, $\ker(\pi_6)$ is a cyclic group and
$\tilde{G}_{\ell}$ could not be its homomorphic image. Let us put
$$G_7=\pi_6(G_6) \subset \PGL_{d'}(k), \quad H_7=\pi_6(H'_6) \subset
G_7.$$ Applying Proposition \ref{cont}(ii)(d2) to $G_6$ and $H_5$,
we conclude that
$$\tilde{G}_{\ell}=G_6/H_6=G_7/H_7.$$
Since $G_7$ is  a subgroup of $\PGL_{d'}(k)\subset
\PGL_{d'}(\bar{\F}_{\ell})$, it follows from a theorem of
Feit-Tits (\cite{FT}; see also \cite{KL}) that $\tilde{G}_{\ell}$
is also isomorphic to a subgroup of $\PGL_{d'}(\bar{\F}_{\ell})$.
In light of property (b),  $d' \ge d$ and we are done.

In particular, if we consider the faithful $\g_j$-module $W_j$
then we get
$$\dim_E(W_j) \ge d=\dim_E(V_E)=\prod_{i\in I}\dim_E(W_i).$$
Since each $W_i$ is a faithful $\g_i$-module and therefore
$\dim_E(W_i)>1$, we conclude that the whole set coincides with
singleton $\{j\}$. This means that $\g_E=g_j$ is absolutely simple
and therefore $\g$ is also absolutely simple. We also conclude
that $V_E=W_j$ has the minimum dimension $d$.

It follows from Weyl's character formula \cite{Bourbaki} that every faithful
simple $\g_E=\g_j$-module of minimum $E$-dimension $d_0$ is
fundamental. This implies that $V_E=W_j$ is fundamental.
\end{proof}

\begin{cor}
\label{Sp} We keep all the notations and assumptions of Theorem
\ref{Lie}. Recall that $\G^0$ is the identity component for
$G_{\alg}$. Assume, in addition that there exists a non-degenerate
alternating bilinear form
$$e: V_{\ell}\times V_{\ell}\to\Q_{\ell}$$
such that $G\subset \Aut(V_{\ell},e)$. We write $\Sp_{V_{\ell},e}$
for the corresponding symplectic $\Q_{\ell}$-algebraic group. If
$d\ne 56$ then $\G=\Sp_{V_{\ell},e}$. If $d=56$ then either
$\G=\Sp_{V_{\ell},e}$ or $\G^0$ is a simply-connected absolutely
simple $\Q_{\ell}$-algebraic group of type $\E_7$.
\end{cor}

\begin{proof}
Clearly, $\G \subset \Sp_{V_{\ell},e}$. This implies that the
$\g$-module $V_{\ell}$ is symplectic and therefore the
$\g_E$-module $V_E$ is also symplectic. By Theorem \ref{Lie},
$\g_E$ is absolutely simple and  $V_E$ is {\sl fundamental} of
minimum dimension. It follows from Tables in \cite{Bourbaki} that
if $V_E$ is a  symplectic fundamental representation of minimum
dimension of $\g_E$ then either $\g_E$ is the Lie algebra of the
symplectic group of $V_E$ or $\dim_E(V_E)=56$ and $\g_E$ is a Lie
algebra of type $\E_7$ and the highest weight of $V_E$ is the only
minuscule dominant weight. One has only to recall that
$$d=\dim_{\Q_{\ell}}(V_{\ell})=\dim_E(V_E).$$
\end{proof}

\begin{cor}
\label{Sp2}
 We keep all the notations and assumptions of Theorem
\ref{Lie}. Assume, in addition that there exists a non-degenerate
alternating bilinear form
$$e: V_{\ell}\times V_{\ell}\to\Q_{\ell}$$
such that $G\subset \Aut(V_{\ell},e)$. We write $\Sp_{V_{\ell},e}$
for the corresponding symplectic $\Q_{\ell}$-algebraic group.
Assume also that $\ell=2$ and $\tilde{G}_{2}$ contains a subgroup
isomorphic to the alternating group $\A_{d+1}$ (e.g.,
$\tilde{G}_{2}\cong \A_{d+1}$ or $\tilde{G}_{2}\cong \A_{d+2}$).
Then $\G=\Sp_{V_{\ell},e}$.
\end{cor}

\begin{proof}
In light of Cor. \ref{Sp} we may assume that $d=56$ and $\G^0$ is
a simply-connected absolutely simple $\Q_{2}$-algebraic group of
type $\E_7$. We have to arrive to a contradiction. First, notice
that in the notations of Step 3 (Proof of Th. \ref{Lie}) $V_E=W_J$
and $\G^0=\G=\prod_{i\in I}\G_i=\G_j$. Also, $\deg(\Ad_{j,E})=2$.
By Step 5 of Proof of Th. \ref{Lie}, there are a compact subgroup
$G_4\subset\G_j(E)=\G^0(E)$, an  open normal subgroup $H_4 \subset
G_4$ such that $G_4/H_4=\tilde{G}_{2}$ and $2$ is the only  prime
divisor of the (super)order of $H_4$.
 Let us put $p=5$. Notice that $p$ is
 {\sl not} a {\sl torsion number} for $\E_7$ (\cite{SerreB},
 1.3.6; \cite{St}).
 Clearly, every Sylow-$p$-subgroup of $G_4$ is is isomorphic to a
Sylow-$p$-subgroup of $\tilde{G}_{2}$. By assumption,
$\tilde{G}_{2}$ contains a subgroup isomorphic to
$\A_{d+1}=\A_{57}$. This implies that $\tilde{G}_{2}$ contains an
elementary abelian $p$-subgroup  $(\Z/p\Z)^8$ and therefore $G_4$
also contains an elementary abelian $p$-subgroup $(\Z/p\Z)^8$.
Since $G_4 \subset \G(E)$, the group $\G^0(\bar{\Q}_{2})$ also
contains $(\Z/p\Z)^8$. Since $p$ is not a torsion number for
$\G^0$, a maximal torus of $\G^0$ contains a subgroup isomorphic
to $(\Z/p\Z)^8$. Since the rank of $\G^0$ is $7<8$, we obtain the
desired contradiction.
\end{proof}

\begin{thm}
\label{Lie2} Suppose  $V_{2}$ is a $\Q_{2}$-vector space of even
finite dimension $d \ge 8$, $G \subset \SL(V_{2})\subset
\Aut(V_{2})$ a compact linear group. Suppose that there exists a
non-degenerate alternating bilinear form
$$e: V_{\ell}\times V_{2}\to\Q_{2}$$
such that $G\subset \Aut(V_{2},e)$.
 We write $G_{alg}$ for the algebraic
envelope of $G$ and $\g$ for the Lie algebra of $G$.

 Suppose $T\subset V_{2}$ is a $G$-stable
$\Z_{2}$-lattice in $V_{2}$. Let us put $V(2):=T/2 T$. Clearly,
$V(2)$ is a $d$-dimensional vector space and carries a natural
structure of $G$-module. Let us denote by $\tilde{G}_{2}$ the
image of the natural homomorphism
$$G \to \Aut_{\Z_{2}}(T) \twoheadrightarrow \Aut(T/2 T)=\Aut_{\F_{2}}(V(2)).$$
Suppose that  $\tilde{G}_{2}$ contains a subgroup isomorphic to
$\A_{d+1}$.

 Then $\G=\Sp_{V_{2},e}$ and $g$ coincides with the Lie algebra
 $\sp_{V_{2},e}$ of the symplectic group $\Sp_{V_{2},e}$.
\end{thm}

\begin{proof}
Clearly, $\G \subset \Sp_{V_{2},e}$. Replacing $G$ by its open
subgroup of finite index, we may assume that
  $\tilde{G}_{2}=\A_{d+1}$.
Taking into account Remark \ref{repA}, we observe that $\ell=2$,
$G$, $d$ and ${G}_{2}$ satisfy the conditions of Theorem \ref{Lie}
and Corollary \ref{Sp2}.
\end{proof}

\section{Applications to abelian varieties}

\begin{thm}
\label{apps} Suppose $K$ is a field with $\fchar(K)\ne 2$, suppose
$X$ is an abelian variety over a $K$ and $\lambda$ is a
polarization on $X$. Suppose $\tilde{G}_{2,X}$ is the image of
$\Gal(K)$ in $\Aut(X_2)$. Let us put $d=2\dim(X)$.  Assume that
$\tilde{G}_{2,X}$ contains a simple non-abelian subgroup
$\mathcal{G}$, enjoying the following properties.

Every faithful projective representation of $\mathcal{G}$ in
characteristic $2$ has dimension $\ge d$. If $d=56$ (i.e.,
$\dim(X)=28$), assume $\mathcal{G}$ contains $\A_{d+1}$. Then:
\begin{enumerate}
\item[(i)]
Either $g_{2,X}=\sp(V_2(X),e_{\lambda})$ or
$g_{2,X}=\Q_2\I\oplus\sp(V_2(X),e_{\lambda})$ or $g_{2,X}=\{0\}$
or $g_{2,X}=\Q_2\I$.

If every finite algebraic extension of $K$ contains only finitely
many $2$-power roots of unity then either
 $g_{2,X}=\Q_2\I$ or $g_{2,X}=\Q_2\I\oplus\sp(V_2(X),e_{\lambda})$.
 If there exists a finite algebraic extension of $K$ that contains all
 $\ell$-power roots of unity then either $g_{2,X}=\sp(V_{2}(X),e_{\lambda})$ or
$g_{2,X}=\{0\}$.
\item[(ii)]
Assume that for each central extension (short exact sequence)
$$1 \to \Z/2\Z \hookrightarrow \mathcal{G}' \twoheadrightarrow
\mathcal{G}\to 1$$ there does not exist such a lifting of the
$\F_{2}[\mathcal{G}']$-module $X_2$ to an absolutely simple
$\Q_{2}[\mathcal{G}']$-module $\Q_{2}^d$ that the distinguished
central subgroup $\Z/2\Z$ acts on $\Q_{2}^d$ via multiplications
by $\pm 1$ and there exists a non-degenerate
$\mathcal{G}'$-invariant alternating bilinear form on $\Q_{2}^d$.
Then:
\begin{enumerate}
\item[(a)]
 the ring $\End(X)$ of all $K_a$-endomorphisms of $X$ is
$\Z$.

\item[(b)]
either $g_{2,X}=\sp(V_2(X),e_{\lambda})$ or
$g_{2,X}=\Q_2\I\oplus\sp(V_2(X),e_{\lambda})$.

 If every finite
algebraic extension of $K$ contains only finitely many $2$-power
roots of unity then
 $g_{2,X}=\Q_2\I\oplus\sp(V_2(X),e_{\lambda})$.
 If there exists a finite algebraic
extension of $K$ that contains all $\ell$-power roots of unity
then $g_{2,X}=\sp(V_{2}(X),e_{\lambda})$.
\end{enumerate}
\end{enumerate}
\end{thm}

\begin{proof}
Replacing if necessary $K$ by its suitable finite separable
extension we may assume that $\tilde{G}_{2,X}=\mathcal{G}$. Let us
put $V_2=V_2(X), T=T_2(X), e=e_{\lambda}$ and
 $$G=\rho_{2,X}(\Gal(K_a/K(2))\subset
 G_{2,X}=\rho_{2,X}(\Gal(K))\subset \Aut_{\Z_2}(T) \subset
 \Aut_{\Q_2}(V_2).$$
  Clearly, $G$ is a normal closed subgroup of $G_{2,X}$
 and the quotient $G_{2,X}/G$ is a one-dimensional or zero-dimensional
 compact commutative $\ell$-adic Lie group.  Since  $\mathcal{G}$
 is simple non-abelian, the image of $G$ in $\Aut(T/2T)$ coincides
 with $\tilde{G}_{2,X}=\mathcal{G}$.

Assume that the condition (ii) holds. Applying Corollaries
\ref{Sp} and \ref{Sp2},
we conclude that the
Lie algebra of $G$ coincides with $\sp(V_2(X),e_{\lambda})$. Since
$G$ is a closed subgroup of $G_{2,X}$,
$$\sp(V_2(X),e_{\lambda})\subset g_{2,X} \subset \Q_2\I\oplus
\sp(V_2(X),e_{\lambda})$$ and therefore either
$g_{2,X}=\sp(V_2(X),e_{\lambda})$ or $g_{2,X}=\Q_2\I\oplus
\sp(V_2(X),e_{\lambda})$. In order to prove that $\End(X)=\Z$,
recall that $\End(X)\otimes \Q_2 \subset \End_{\Q_2}V_2(X)$
commutes with $g_{2,X}$. Since the centralizer of
$\sp(V_2(X),e_{\lambda})$ consists of scalars, $\End(X)\otimes
\Q_2=\Q_2\I$ and therefore $\End(X)=\Z$.

Assume that the condition (i) holds. If $G$ is infinite then
applying Theorem \ref{Lie2}, we conclude by the same token that
either $g_{2,X}=\sp(V_2(X),e_{\lambda})$ or $g_{2,X}=\Q_2\I\oplus
\sp(V_2(X),e_{\lambda})$. Assume that $G$ is finite. Then the
dimension of $G_{2,X}$ is either $0$ or $1$. If it is $0$ then
$g_{2,X}=\{0\}$. Assume that the dimension of $G_{2,X}$ is $1$ and
therefore $G_{2,X}$ is infinite and $g_{2,X}$ is a one-dimensional
$\Q_2$-vector subspace of $\End_{\Q_2}(V_2(X)$.

Our goal now is to prove that $g_{2,X}$ consists of scalars. Let
$u$ be a non-zero element of $g_{2,X}$. Clearly, $g_{2,X}=\Q_2 u$
and the conjugation by $G_{2,X}$ leaves $g_{2,X}=\Q_2 u$ stable.
This implies that the conjugation by $G_{2,X}$ also leaves stable
the subalgebra $\Q_2[u]\subset \End_{\Q_2}(V_2(X))$. Assume for a
moment that the $G_{2,X}$-module $V_2(X)$ is very simple. Then
either $\Q_2[u]= \End_{\Q_2}(V_2(X))$ or $\Q_2[u]=\Q_2\I$. The
former equality could not be true, since $\Q_2[u]$ is commutative
while $\End_{\Q_2}(V_2(X))$ is not. Hence $\Q_2[u]=\Q_2\I$, i.e.,
$u$ is a scalar and therefore $g_{2,X}=\Q_2 u=\Q_2\I$ and we are
done.

In order to finish the proof, recall that if every finite
algebraic extension of $K$ contains only finitely many $2$-power
roots of unity then $g_{2,X}$ does {\sl not} lie in $\sll(V_2(X))$
and therefore $g_{2,X}$ is neither $\{0\}$
 nor $\sp(V_2(X),e_{\lambda})$.
 If there exists a finite algebraic extension of $K$ that contains all
 $\ell$-power roots of unity then
$$g_{2,X} \subset \sll(V_2(X))\bigcap
[\Q_2\I\oplus\sp(V_2(X),e_{\lambda})]=\sp(V_2(X),e_{\lambda})$$
and therefore $g_{2,X}$ is neither
$\Q_2\I\oplus\sp(V_2(X),e_{\lambda})$ nor $\Q_2\I$.
\end{proof}

\begin{cor}
\label{appsM} Suppose $K$ is a field with $\fchar(K)\ne 2$,
suppose $X$ is a $5$-dimensional abelian variety over a $K$ and
$\lambda$ is a polarization on $X$. Suppose $\tilde{G}_{2,X}$ is
the image of $\Gal(K)$ in $\Aut(X_2)$.  Assume that
$\tilde{G}_{2,X}$ contains a subgroup isomorphic to $\M_{11}$
(e.g., $\tilde{G}_{2,X}=\M_{11}$ or $\tilde{G}_{2,X}$ is
isomorphic to the Mathieu group $\M_{12}$). Then:
\begin{enumerate}
\item[(a)]
 the ring $\End(X)$ of all $K_a$-endomorphisms of $X$ is
$\Z$.

\item[(b)]
either $g_{2,X}=\sp(V_2(X),e_{\lambda})$ or
$g_{2,X}=\Q_2\I\oplus\sp(V_2(X),e_{\lambda})$.

 If every finite
algebraic extension of $K$ contains only finitely many $2$-power
roots of unity then
 $g_{2,X}=\Q_2\I\oplus\sp(V_2(X),e_{\lambda})$.
 If there exists a finite algebraic
extension of $K$ that contains all $\ell$-power roots of unity
then $g_{2,X}=\sp(V_{2}(X),e_{\lambda})$.
\end{enumerate}
\end{cor}

\begin{proof}
We have $d:=2\dim(X)=10$. The proof follows readily from Theorem
\ref{apps}(ii) (applied to $\mathcal{G}=\M_{11}$) combined with
Remark \ref{repM}.
\end{proof}

\begin{cor}
\label{appsL} Suppose $K$ is a field with $\fchar(K)\ne 2$,
suppose $X$ is an abelian variety over a $K$ and $\lambda$ is a
polarization on $X$. Suppose $\tilde{G}_{2,X}$ is the image of
$\Gal(K)$ in $\Aut(X_2)$. Let us put $d=2\dim(X)$.  Assume that
$\tilde{G}_{2,X}$ contains a simple non-abelian subgroup
$\mathcal{G}$, enjoying the following properties.

There exist an odd power prime $q$ and an integer $m \ge 3$ such
that $(q,m) \ne (3,4)$,
 $\frac{q^m-1}{q-1}=d+1$ or $d+2$ and $\mathcal{G} \cong
 \L_m(q)$. Then:
\begin{enumerate}
\item[(a)]
 the ring $\End(X)$ of all $K_a$-endomorphisms of $X$ is
$\Z$.

\item[(b)]
either $g_{2,X}=\sp(V_2(X),e_{\lambda})$ or
$g_{2,X}=\Q_2\I\oplus\sp(V_2(X),e_{\lambda})$.

 If every finite
algebraic extension of $K$ contains only finitely many $2$-power
roots of unity then
 $g_{2,X}=\Q_2\I\oplus\sp(V_2(X),e_{\lambda})$.
 If there exists a finite algebraic
extension of $K$ that contains all $\ell$-power roots of unity
then $g_{2,X}=\sp(V_{2}(X),e_{\lambda})$.
\end{enumerate}
\end{cor}

\begin{proof}
The proof follows readily from Theorem \ref{apps}(ii) (applied to
$\mathcal{G}=\L_m(q)$) combined with Remark \ref{repL}.
\end{proof}

\section{Proof of main results}
 \label{mainproof}

\begin{proof}[Proof of Theorem \ref{main0}]
In order to get the assertions about $g_{2,X}$, one has only to
combine Remark \ref{repA} and Theorem \ref{apps} (applied to
$\mathcal{G}=\A_{d+1}$).
\end{proof}

\begin{rem}
\label{jac2} Suppose $f(x)\in K[x]$ is a polynomial of degree $n
\ge 5$ without multiple roots and $X=J(C_f)$ is the jacobian of
$C=C_f: y^2=f(x)$. One may easily check that $d=2\dim(J(C))=n-1$
if $n$ is odd and $d=2\dim(J(C))=n-2$ when $n$ is even. It is
well-known (see for instance Sect. 5 of \cite{ZarhinTexel}) that
if $X=J(C_f)$ is the jacobian of
$$C=C_f: y^2=f(x)$$
where $f(x) \in K[x]$ has no multiple roots then $\tilde{G}_{2,X}$
is isomorphic to $\Gal(f)$. Clearly, if $n\ge 9$  and $\Gal(f)$
contains $\A_n$ then $\Gal(f)$ contains $\A_{d+1}$. It is also
clear that if $n=12$ and $\Gal(f)$ contains $\M_{12}$ then it also
contains $\M_{11}$.

\end{rem}
\begin{proof}[Proof of Theorem \ref{main}]
It follows from Remark \ref{jac2} combined with Theorem \ref{apps}
and Corollaries  \ref{appsM}  and  \ref{appsL} that $\End(X)=\Z$
and $g_{2,X}$ contains $\sp(V_2(X),e_{\lambda})$. One has only to
recall that $\sp(V_2(X),e_{\lambda})$ coincides with the Lie
algebra of $\Aut_{\Z_2}(T_2(X), e_{\lambda})$ and $g_{2,X}$ is the
Lie algebra of $\rho_{2,J(C_f)}(\Gal(K))$.
\end{proof}

\begin{proof}[Proof of Theorem \ref{main1}]
 By Th. \ref{main0},
$g_{2,X}=\Q_2\I\oplus\sp(V_2(X),e_{\lambda})$.
We need  to prove that
$g_{\ell,X}=\Q_{\ell}\I\oplus\sp(V_{\ell}(X),e_{\lambda})$ for all
primes $\ell \ne \fchar(K)$. This assertion follows readily from
the next auxiliary statement.

\begin{lem}
\label{rank} Assume that the field $K$ is either finitely
generated over $\Q$ or a global field of characteristic $>2$. Let
$X$ be an abelian variety defined over $K$  such that
$g_{2,X}=\Q_2\I\oplus\sp(V_2(X),e_{\lambda})$. Then
$g_{\ell,X}=\Q_{\ell}\I\oplus\sp(V_{\ell}(X),e_{\lambda})$ for all
primes $\ell \ne \fchar(K)$.
\end{lem}

\begin{proof}[Proof of Lemma \ref{rank}]

First, assume that $K$ is a global field. It follows easily (as in
the proof of Th. \ref{apps}(ii)) that $\End(X)=\Z$. Now it follows
from results of \cite{ZarhinTor}, \cite{Faltings} that
$g_{\ell,X}$ is an absolutely irreducible reductive subalgebra of
$\End_{\Q_{\ell}}(V_{\ell}(X))$ for all primes $\ell \ne
\fchar(K)$. We also have
$$g_{\ell,X}\subset \Q_{\ell}\I\oplus\sp(V_{\ell}(X),e_{\lambda})\subset \End_{\Q_{\ell}}(V_{\ell}(X)).$$
Notice that the rank of the reductive $\Q_{\ell}$-Lie algebra
$g_{\ell,X}$ does not depend on the choice of $\ell$
\cite{ZarhinInv}. This implies that $g_{\ell,X}$ and
$\Q_{\ell}\I\oplus\sp(V_{\ell}(X),e_{\lambda})$ have the same rank
$\dim(X)+1$. It follows from a variant of a theorem of Borel - de
Siebenthal (\cite{ZarhinDuke}, Key Lemma on p. 522) that
$g_{\ell,X}= \Q_{\ell}\I\oplus\sp(V_{\ell}(X),e_{\lambda})$. This
proves Lemma \ref{rank} in the case of global $K$.

The case of of an arbitrary $K$ finitely generated over $\Q$
follows  from the case of number field with the help of Serre's
variant of Hilbert irreducibility theorem for infinite extensions
in the case of characteristic zero (\cite{SerreMW}, Sect. 10.6;
\cite{SerreRibet}, Sect. 1; \cite{Noot}, Prop. 1.3). Indeed, let
us fix an odd prime $\ell$. Then there exists a number field $K_0$
and an abelian variety $Y$ over $K_0$ such that
$\dim(Y)=\dim(X):=g$ and
$$g_{\ell,X} \cong g_{\ell,Y}, \quad g_{2,X} \cong g_{2,Y}$$
as $\Q_{\ell}$-Lie algebra and $\Q_2$-Lie algebra respectively.
This implies that
$$g_{2,Y} \cong g_{2,X}=\Q_2\I\oplus\sp(V_2(X),e_{\lambda})\cong
\Q_2\I\oplus\sp(2g,\Q_2)$$
 and easy dimension arguments imply that
$Y$ over $K_0$ satisfies the conditions of Lemma \ref{rank}. Since
we already know that Lemma \ref{rank} is true in the case of
ground number field, we conclude that
$$g_{\ell,Y} \cong
\Q_{\ell}\I\oplus\sp(2g,\Q_{\ell})$$
 and therefore
 $$g_{\ell,X}\cong \Q_{\ell}\I\oplus\sp(2g,\Q_{\ell}).$$
Again easy dimension arguments imply that
$$g_{\ell,X}=\Q_{\ell}\I\oplus\sp(V_{\ell}(X),e_{\lambda}).$$
\end{proof}
This ends the proof  of Theorem \ref{main1}.
\end{proof}

\begin{proof}[Proof of Theorem \ref{jac}]
The proof is a straightforward application of Theorem \ref{main1}
combined with Remark \ref{jac2}.
\end{proof}

\begin{proof}[Proof of Theorem \ref{jacall}]
In light of Theorem \ref{jac}, we may assume that $5\le n\le 8$.
Let us put $X=J(C_f)$. Clearly, $\dim(X)=2$ or $3$. Notice also
that $\End(X)=\Z$ (\cite{ZarhinMRL}). Let us put $g:=\dim(X)$.

 First assume that $K$ is a number field. Then it follows from results
of \cite{ZarhinIzv} (combined with the results of \cite{Faltings})
that
$$g_{\ell,X}= \Q_{\ell}\I\oplus\sp(V_{\ell}(X),e_{\lambda})\cong \Q_{\ell}\I\oplus\sp(2g,\Q_{\ell})$$
for all primes $\ell$.

Now assume that $K$ is an arbitrary field of characteristic zero
finitely generated over $\Q$. Let us fix a prime $\ell$. Thanks to
Prop. 1.3 and Cor. 1.5 in \cite{Noot}, there exists a number field
$K_0$ and an abelian variety $Y$ over $K_0$ such that
$\dim(Y)=\dim(X)=g$, $\End(X)\cong \End(Y)$ and $g_{\ell,X} \cong
g_{\ell,Y}$. Since $\End(X)=\Z$, we conclude that $\End(Y)=\Z$.

Since we already know that Theorem \ref{jacall} is true in the
case of ground number field, we conclude that
$$g_{\ell,Y} \cong
\Q_{\ell}\I\oplus\sp(2g,\Q_{\ell})$$
 and therefore
 $$g_{\ell,X}\cong \Q_{\ell}\I\oplus\sp(2g,\Q_{\ell}).$$
 Now easy dimension arguments imply that
$$g_{\ell,X}=\Q_{\ell}\I\oplus\sp(V_{\ell}(X),e_{\lambda}).$$
\end{proof}

\begin{rem}
Concerning specializations of the endomorphism rings of abelian
varieties see also \cite{Masser}.
\end{rem}

\section{Tate classes}
\label{HT}
\begin{thm}
\label{Tate} Let $K$ be a field with $\fchar(K) \ne 2$,
 $K_s$ its separable algebraic closure,
$f(x) \in K[x]$ a separable polynomial of  degree $n \ge 5$, whose
Galois group $\Gal(f)$ enjoys one of the following properties:
\begin{enumerate}
\item[(i)]
$\Gal(f)$ is either $\Sn$ or $\An$;
\item[(ii)]
$n=11$ and $\Gal(f)$ is the Mathieu group $\M_{11}$;
\item[(iii)]
$n=12$ and $\Gal(f)$ is either the Mathieu group $\M_{12}$ or
$\M_{11}$.
\item[(iv)]
There exist an odd power prime $q$ and an integer $m \ge 3$ such
that $(q,m) \ne (3,4)$,
 $n=\frac{q^m-1}{q-1}$ and $\Gal(f)$
contains a subgroup isomorphic  to the projective special linear
group $\L_m(q):=\PSL_m(\F_q)$. (E.g., $\Gal(f)$ is isomorphic
either to the projective linear group $\PGL_m(\F_q)$ or to
$\L_m(q)$.)
\end{enumerate}
 Let $C_f$ be the hyperelliptic curve $y^2=f(x)$ and let
$J(C_f)$ be its jacobian. Assume, in addition, that either $K$ is
a global field of positive characteristic and $n \ge 9$ or $K$ is
a field of characteristic zero finitely generated over $\Q$.

Let $K'\subset K_s$ be a separable finite algebraic extension of
$K$. Then for all primes $\ell \ne \fchar(K)$ and on  each
self-product $J(C_f)^M$ of $J(C_f)$ every $\ell$-adic Tate class
over $K'$ can be presented as a linear combination of products of
divisor classes over $K_a$. In particular, the Tate conjecture is
valid for all $J(C_f)^M$ over all $K'$.
\end{thm}

\begin{proof}
Let us put $X:=J(C_f)$. Recall that one may view $\ell$-adic Tate
classes on  self-products  of $X$ as tensor invariants of
$g_{\ell,X}\bigcap \sp(V_{\ell}(X),e_{\lambda})$ \cite{Tate}.

 Assume that $g_{\ell,X}=\Q_{\ell}\I\oplus\sp(V_{\ell}(X),e_{\lambda})$.
 Then
 $$g_{\ell,X}\bigcap
 \sp(V_{\ell}(X),e_{\lambda})=\sp(V_{\ell}(X),e_{\lambda}).$$

 It follows with the help of results from the invariant theory for symplectic
groups (\cite{Ribet}, \cite{ZarhinEssen}) that each   $\ell$-adic
Tate class on $X^M=J(C_f)^M$ could be presented as a linear
combination of  products of  divisor classes and therefore is
algebraic.
\end{proof}

\begin{rem}
The Tate conjecture is true in codimension $1$ for arbitrary
abelian varieties over $K$ \cite{ZarhinT}, \cite{Faltings2}.
\end{rem}

\section{Hodge classes}
\begin{thm}
\label{Hodge} Suppose $f(x)\in \C[x]$ is a polynomial of  degree
$n \ge 5$ and without multiple roots. Let $C_f$ be the
hyperelliptic curve $y^2=f(x)$ and $J(C_f)$ its jacobian. Assume
that all the coefficients of $f$ lie in a subfield $K \subset \C$
and  the Galois group $\Gal(f)$ of $f$ over $K$ enjoys one of the
following properties:
\begin{enumerate}
\item[(i)]
$\Gal(f)$ is either $\Sn$ or $\An$;
\item[(ii)]
$n=11$ and $\Gal(f)$ is the Mathieu group $\M_{11}$;
\item[(iii)]
$n=12$ and $\Gal(f)$ is either the Mathieu group $\M_{12}$ or
$\M_{11}$.
\item[(iv)]
There exist an odd power prime $q$ and an integer $m \ge 3$ such
that $(q,m) \ne (3,4)$,
 $n=\frac{q^m-1}{q-1}$ and $\Gal(f)$
contains a subgroup isomorphic  to the projective special linear
group $\L_m(q):=\PSL_m(\F_q)$. (E.g., $\Gal(f)$ is isomorphic
either to the projective linear group $\PGL_m(\F_q)$ or to
$\L_m(q)$.)
\end{enumerate}
 Then
each Hodge class on every self-product $J(C_f)^M$ of $J(C_f)$ can
be presented as a linear combination of products of divisor
classes. In particular, the Hodge conjecture is valid for all
$J(C_f)^M$.
\end{thm}
\begin{proof}
If $n\le 8$ then $\dim(J(C_f)) \le 3$. But it is well-known that
the assertion of the theorem is true for all complex abelian
varieties, whose dimension does not exceed $3$ (see, for instance,
\cite{MZ}). Further, we assume that
$$n \ge 9.$$
 Replacing $K$ by its subfield obtained by adjoining to $\Q$ all
coefficients of $f$, we may assume that $K \subset \C$ is finitely
generated over $\Q$. Let us put $X=J(C_f)$. Thanks to Theorem
\ref{jac},
$$g_{\ell,X}=
\Q_{\ell}\I\oplus\sp(V_{\ell}(X),e_{\lambda})$$ for all primes
$\ell$.

The first rational homology group $\Pi_{\Q}:=\H_1(X(\C),\Q)$ of
the complex torus $X(\C)$ is a $2\dim(X)$-dimensional $\Q$-vector
space  provided with a natural structure of rational polarized
Hodge structure of weight $-1$. We refer to \cite{Ribet} for the
definition of its Mumford-Tate group
$\MT=\MT_X\subset\GL(\Pi_{\Q})$. It is a reductive algebraic
$\Q$-group. We write $\mt=\mt_X$ for the Lie algebra of $\MT$; it
is a reductive completely reducible $\Q$-Lie subalgebra of
$\End_{\Q}(\Pi_{\Q})$. The polarization $\lambda$ on $X$ gives
rise to a non-degenerate alternating bilinear (Riemann) form
\cite{MumfordAV}

$$L_{\lambda}: \Pi_{\Q}\times \Pi_{\Q} \to \Q$$ such that
$$\mt \subset \Q\I \oplus\sp(\Pi_{\Q},L_{\lambda}).$$

 Let us choose by $K_a$ the algebraic closure of $K$
in $\C$. For each prime $\ell$ let us put
$$\Pi_{\ell}:=\Pi_{\Q}\otimes_{\Q}\Q_{\ell}.$$
Then there is the well-known natural isomorphism \cite{MumfordAV},
\cite{ZarhinDuke}
 $$\gamma_{\ell}:\Pi_{\ell}\cong V_{\ell}(X)$$ such
that, by a theorem of Piatetski--Shapiro - Deligne - Borovoi
\cite{Deligne}, \cite{SerreK},
$$\gamma_{\ell}g_{\ell,X} \gamma_{\ell}^{-1} \subset
\mt\otimes_{\Q}\Q_{\ell} \subset [\Q\I
\oplus\sp(\Pi_{\Q},L_{\lambda})]\otimes_{\Q}\Q_{\ell} \subset
\End_{\Q_{\ell}}(\Pi_{\ell}).$$

Now easy dimension arguments imply that
$$\gamma_{\ell}g_{\ell,X} \gamma_{\ell}^{-1}=
\mt\otimes_{\Q}\Q_{\ell}= [\Q\I
\oplus\sp(\Pi_{\Q},L_{\lambda})]\otimes_{\Q}\Q_{\ell}$$
 and therefore
$$\mt=\Q\I \oplus\sp(\Pi_{\Q},L_{\lambda}).$$

Recall that Hodge classes on self-products on $X$ could be viewed
as tensor invariants of
$\mt\bigcap\sp(\Pi_{\Q},L_{\lambda})=\sp(\Pi_{\Q},L_{\lambda})$.
As in the case of Tate classes, results from the invariant theory
for symplectic groups (\cite{Ribet} imply that each Hodge class on
a self-product of $X$ can be presented a linear combination of
products of divisor classes.
\end{proof}

\begin{rem}
In the course of the proof we established the equality
$\gamma_{\ell}g_{\ell,X} \gamma_{\ell}^{-1}=
\mt\otimes_{\Q}\Q_{\ell}$. This proves the Mumford-Tate conjecture
\cite{SerreK}  for $X=J(C_f)$ in the case of finitely generated
$K$ when $n \ge 9$.
\end{rem}

\bigskip

\noindent {\small {Department of Mathematics, Pennsylvania State University,}}

\noindent {\small {University Park, PA 16802, USA} }

\vskip .4cm

\noindent {\small {Institute for Mathematical Problems in Biology,}}

\noindent {\small {Russian Academy of Sciences, Push\-chino, Moscow Region, 142292, RUSSIA}}

 \vskip  .4cm

\noindent {\small {\em E-mail address}: zarhin@math.psu.edu}

\end{document}